\documentclass[aos,preprint]{imsart}
\usepackage{paralist}
\usepackage{enumerate} 
\usepackage{dcolumn}
\usepackage{floatrow}
\usepackage[OT1]{fontenc}
\usepackage[numbers]{natbib}
\usepackage[colorlinks,citecolor=blue,urlcolor=blue]{hyperref}
\usepackage{amsmath,amssymb,amsthm,amscd,graphicx,color,accents}

\usepackage{caption}
\usepackage{lipsum}
\usepackage{floatrow}
\usepackage{lipsum}
\usepackage{placeins}

\startlocaldefs
\newtheorem{theorem}{Theorem}[section]
\newtheorem{proposition}[theorem]{Proposition}
\newtheorem{lemma}[theorem]{Lemma}
\newtheorem{corollary}[theorem]{Corollary}

\newtheorem{definition}[theorem]{Definition}
\newtheorem{example}[theorem]{Example}

\newcommand{\Pro}{\par\noindent\rm\textsc{Proof. }}

\def\tr{{\rm{trace}\,}}

\def\bX{{\boldsymbol{X}}}
\def\bY{{\boldsymbol{Y}}}
\def\bZ{{\boldsymbol{Z}}}
\def\bx{{\boldsymbol{x}}}
\def\by{{\boldsymbol{y}}}

\def\dd{{{\hskip1pt}\rm{d}}}

\def\P{{\mathbb P}}

\def\U{{\mathcal U}}
\def\V{{\mathcal V}}
\def\E{\mathbb{E}}

\def\C{\mathbb{C}}
\def\N{\mathbb{N}}
\def\R{\mathbb{R}}

\def\Cov{{\rm{Cov}}}
\def\Cor{{\rm{Cor}}}
\def\Var{{\rm{Var}}}

\def\cid{\stackrel{\mbox{\tiny d}}{\longrightarrow}}
\def\cip{\stackrel{\mbox{\tiny P}}{\longrightarrow}}
\numberwithin{equation}{section}
\theoremstyle{plain}

\endlocaldefs

\begin{document}
	\begin{frontmatter}
		
		\title{The Distance Standard Deviation}
		\runtitle{The Distance Standard Deviation}
		
		\begin{aug}
			\author{\fnms{Dominic} \snm{Edelmann}\thanksref{t1}\ead[label=e1]{dominic.edelmann@dkfz-heidelberg.de},}
			\author{\fnms{Donald} \snm{Richards}\ead[label=e2]{richards@stat.psu.edu},}
			\and
			\author{\fnms{Daniel} \snm{Vogel}\ead[label=e3]{daniel.vogel@abdn.ac.uk}}
			
			\thankstext{t1}{Corresponding author.}
			\runauthor{D. Edelmann, D. Richards, and D. Vogel}
			
			\affiliation{German Cancer Research Center, Pennsylvania State University, and University of Aberdeen}
			
			\address{German Cancer Research Center\\ %(Deutsches Krebsforschungszentrum), 
				Im Neuenheimer Feld 280\\
				69120 Heidelberg, Germany\\
				\printead{e1}\\}
			
			\address{Department of Statistics\\ 
				Pennsylvania State University\\ 
				University Park, PA 16802, U.S.A.\\
				\printead{e2}\\}
			
			\address{Institute for Complex Systems and Mathematical Biology\\ 
				University of Aberdeen\\ 
				Aberdeen AB24 3UE, U.K.\\
				\printead{e3}
			}			
		\end{aug}
		
		\begin{abstract}
			The distance standard deviation, which arises in distance correlation analysis of multivariate data, is studied as a measure of spread. The asymptotic distribution of the empirical distance standard deviation is derived under the assumption of finite second moments. Applications are provided to hypothesis testing on a data set from materials science and to multivariate statistical quality control. The distance standard deviation is compared to classical scale measures for inference on the spread of heavy-tailed distributions.  
				Inequalities for the distance variance are derived, proving that the distance standard deviation is bounded above by the classical standard deviation and by Gini's mean difference. New expressions for the distance standard deviation are obtained in terms of Gini's mean difference and the moments of spacings of order statistics. It is also shown that the distance standard deviation satisfies the axiomatic properties of a measure of spread.  
		\end{abstract}

		
		
	\end{frontmatter}

	\section{Introduction}
	\label{sec:intro}
	
	In recent years, the topic of distance correlation has been prominent in statistical analyses of dependence between multivariate data sets.   The concept of distance correlation was defined by Sz\'ekely, Rizzo, and Bakirov \citep{szekely2007} and Sz\'ekely and Rizzo \citep{szekely2009}, and they applied distance correlation methods to testing independence and measuring association between collections of random vectors.
	
	Since the appearance of the papers \citep{szekely2007,szekely2009}, enormous interest in the theory and applications of distance correlation has arisen.  We refer to the articles \citep{rizzo2010,%szekely2012,
		szekely2013,szekely2014} on statistical inference, \citep{,
		fiedler2016a,fokianos2016,jentsch2016,zhou2012} on time series, \citep{Dueck2014,Dueck2015,Dueck2016} on affinely invariant distance correlation and connections with singular integrals, \citep{lyons2013} on metric spaces, and \citep{sejdinovic2013} 
	on machine learning.  Distance correlation methods have also been applied 
	to assessing associations between familial relationships, lifestyle factors, diseases, and mortality \citep{kong2012}, and to detecting associations in large astrophysical databases \citep{martinez2014,richards2014}.

	For $z \in \C$, denote by $|z|$ the modulus of $z$.  For a positive integer $p$ and $s,x \in \R^p$, denote by $\langle s,x\rangle$ the Euclidean inner product on $\R^p$ and by $\|s\| = \langle s,s \rangle^{1/2}$ the corresponding Euclidean norm.  We also define the constant 
	$$
	c_p = \frac{\pi^{(p+1)/2}}{\Gamma\big((p+1)/2\big)}.
	$$
For random vectors $X \in \R^p$ and $Y \in \R^q$, let 
	$$
	f_{X,Y}(s,t) = \E \exp\big(\sqrt{-1}(\langle s,X\rangle + \langle t,Y\rangle) \big),
	$$
where	$s \in \R^p$, $t \in \R^q$, be the joint characteristic function of $(X,Y)$ and let $f_X(s) = f_{X,Y}(s,0)$ and $f_Y(t) = f_{X,Y}(0,t)$ be the corresponding marginal characteristic functions.  The {\it distance covariance} between $X$ and $Y$ is defined as the nonnegative square root of 
	\begin{equation}
	\label{eq:dcov}
	\V^2(X,Y) = \frac{1}{c_p c_q} \int_{\R^{p+q}} 
	\big|f_{X,Y}(s,t)-f_X(s)f_Y(t)\big|^2 \, \frac{\dd s {\hskip 1pt}\dd t}{\|s\|^{p+1} \, \|t\|^{q+1}},
	\end{equation}
	the {\it distance variance} is defined as
	\begin{align} \label{eq:dvar}
	\V^2(X) = \V^2(X,X) &= \frac{1}{c_p^2} \int_{\R^{2p}} 
	\big|f_{X}(s+t)-f_X(s)f_X(t)\big|^2 \, \frac{\dd s {\hskip 1pt} \dd t}{\|s\|^{p+1} \, \|t\|^{p+1}},
	\end{align}
	and the {\it distance standard deviation}, $\V(X)$, is defined as the nonnegative square root of $\V^2(X)$. (We note that this terminology differs from that of Sz\'{e}kely, et al. \citep{szekely2009,szekely2007}, who refer to $\V(X)$ as the {\it distance variance}; we will refer to $\V(X)$ instead as the {\it distance standard deviation}, which is justified by the fact that $\V(X)$ satisfies an equivariance property that is given below in \eqref{eq:dilation}.)
Also, the {\it distance correlation coefficient} is defined as
	\begin{equation}
	\label{eq:dcor}
	\mathcal{R}(X,Y) = \frac{\V(X,Y)}{\sqrt{\V(X) \V(Y)}}
	\end{equation}
	as long as $\V(X), \V(Y) \neq 0$, and zero otherwise.  We remark that the weighted $L^2$-norm in (\ref{eq:dcov}) was studied in the univariate setting by Feuerverger \citep{feuerverger1993}. 
	
The distance correlation coefficient, unlike the Pearson correlation coefficient, characterizes independence: $\mathcal{R}(X,Y)= 0$ if and only if $X$ and $Y$ are mutually independent.  Moreover, $0 \leq \mathcal{R}(X,Y) \leq 1$ and, for one-dimensional random variables $X,Y \in \R$, we have $\mathcal{R}(X,Y)=1$ if and only if $Y$ is a linear function of $X$, almost surely.  The empirical distance correlation possesses a remarkably simple expression \citep[Theorem 1]{szekely2007}, and efficient algorithms for computing it are now available \citep{Huo2015}.

We note that $\mathcal{R}(X,Y)$ is one of several coefficients characterizing independence that are applicable to hypothesis testing.  Other concepts of dependence are, e.g., the Hilbert-Schmidt Independence Criterion (HSIC) \cite{gretton2008}, ball covariance \cite{pan2019} and mutual information \cite{Berret2019}. Each of these concepts satisfy numerous desirable properties, and the comparison of their properties and finite-sample performance is an active area of research \cite{ramdas2015decreasing,sejdinovic2013,simon2014comment}.

An interesting property of the distance covariance is that its square is translation-invariant and scale-equivariant, which implies that the distance standard deviation satisfies
	\begin{equation} \label{eq:dilation}
	\V(a+bX) = |b| \, \V(X),
	\end{equation}
for all $a,b \in \R$ \citep[Theorem 4]{szekely2007}. Moreover, $\V(X)$ is defined for all random variables $X$ with finite first moments, whereas the classical standard deviation requires the existence of finite second moments. These properties suggest that the distance standard deviation is a potentially interesting measure of scale for heavy-tailed distributions. As the term $\V(X) \, \V(Y)$ appears in the denominator of $\mathcal{R}(X,Y)$, a study of properties of the distance standard deviation may lead to a better understanding of the distance correlation.

In this paper, we study the distance standard deviation $\V(X)$ and provide applications to hypothesis testing and multivariate statistical quality control. We apply the distance standard deviation to a data set, originating from materials science, on a physical model for describing a crystal undergoing a structural phase transition when subjected to several cooling-heating cycles.
In a different direction, we further show how $\V(X)$ can be applied in the statistical quality control of multivariate production processes.

We will also compare $\V(X)$ to other measures of spread.  Indeed, suppose that $\E (\|X\|^2) < \infty$, and let $X$, $X'$, and $X''$ be independent and identically distributed (i.i.d.); then, by Sz\'ekely, et al. \citep[Remark 3]{szekely2007}, 
\begin{equation} 
\label{rep:dvar2}
\V^2(X) = \E (\|X-X'\|^2) + (\E \|X-X'\|)^2 - 2 \E (\|X-X'\| \cdot \|X-X''\|).
\end{equation}
The second term on the right-hand side of (\ref{rep:dvar2}) is reminiscent of Gini's mean difference \citep{gerstenberger2015,yitzhaki2013}, which is defined for real-valued random variables $Y$ as
\begin{equation} \label{GMD}
\Delta(Y) =	\E |Y-Y'|,
\end{equation}
where $Y$ and $Y'$ are i.i.d.  Furthermore, if $X \in \R$ then one-half the first summand in	(\ref{rep:dvar2}) equals $\sigma^2(X)= E(X^2) - E(X)^2$, the variance of $X$.

We provide a detailed comparison of $\V(X)$, $\Delta(X)$, and $\sigma(X)$.  We demonstrate that when the distributions of interest are heavy-tailed, $\V(X)$ provides estimators of scale that are asymptotically more efficient than estimators based on $\Delta(X)$ or $\sigma(X)$.  Moreover, several inequalities between $\V(X)$, $\Delta(X)$, and $\sigma(X)$ are derived.

We further show that the distance standard deviation is an axiomatic measure of spread in the sense of Bickel and Lehmann \citep{Bickel2012}. According to 
{\citep{Bickel2012}}, a {\it measure of spread} is a functional $\tau(X)$ satisfying the axioms: 
\begin{enumerate}
\setlength{\itemsep}{0pt}
\item[(C1)] $\tau(X) \geq 0$,
\item[(C2)] $\tau(a+bX) = |b| \, \tau (X)$ for all $a,b \in \R$, and 
\item[(C3)] $\tau(X) \leq \tau(Y)$ if for all $0 < \alpha \leq \beta < 1$,
\begin{equation*}
%\label{eq:dispersive}
	F^{-1}(\beta) - F^{-1}(\alpha) \leq  G^{-1}(\beta) - G^{-1}(\alpha),
	\end{equation*}
	where $F$ and $G$ are the cumulative distribution functions of $X$ and $Y$, respectively, and $F^{-1}$ and $G^{-1}$ are the corresponding right-continuous inverses.
\end{enumerate}
The distance covariance obviously satisfies (C1) and (C2). We will show that $\V(X)$ also satisfies (C3), hence proving that $\V(X)$ is a measure of spread in the above sense.  However, we will also establish some clear differences between $\V(X)$, on the one hand, and $\Delta(X)$ and $\sigma(X)$, on the other hand.

	The paper is organized as follows. In Section \ref{sec:estimator}, the asymptotic distribution of the empirical distance standard deviation under the existence of the second moment of $X$ is derived. The asymptotic relative efficiency (ARE) of the empirical distance standard deviation with respect to competing estimators of spread is evaluated for various distributions. 
	In Section \ref{sec:app}, we apply the empirical distance standard deviation to perform two-sample hypothesis testing for a data set from materials science and we also show the applicability of the empirical distance standard deviation in multivariate statistical quality control.  Further, we demonstrate the superior performance of tests based on the distance standard deviation when the underlying distributions are heavy-tailed. 		
	In Section \ref{sec:ineq}, we derive inequalities between the summands in the distance variance representation (\ref{rep:dvar2}).  We will prove in the case of scalar random variables that $\V(X)$ is bounded above by $\Delta(X)$ and by $\sigma(X)$. In Section \ref{sec:properties}, we show that the representation (\ref{rep:dvar2}) can be simplified further, revealing relationships between $\V(X)$ and the moments of spacings of order statistics. Using novel representations, we show that $\V(X)$ is a measure of spread in the sense of \citep{Bickel2012}; moreover, we identify crucial differences between $\V(X)$, $\Delta(X)$ and $\sigma(X)$. 
	% Section \ref{sec:distr} provides closed-form expressions for the distance variance for numerous parametric distributions. 
	We conclude the paper in Section \ref{sec:discussion} with a discussion of the given results. All proofs are provided in the supplementary material.

	\section{The empirical distance standard deviation}
	\label{sec:estimator}

	In order to develop an empirical version of $\V^2(X)$, Sz\'{e}kely, et al. \citep{szekely2007,szekely2009} derived an alternative representation of $\V^2(X)$; they showed that 
%		The integral representation in equation (\ref{eq:dvar}) of the distance variance $\V^2(X)$ generally is not suitable for practical purposes.  Sz\'{e}kely, et al. \citep{szekely2007,szekely2009} derived an alternative representation; they show that 
if the random vector $X \in \R^p$ satisfies $\E \|X\|^2 < \infty$ and if $X$, $X'$, and $X''$ are i.i.d. then 
	\begin{equation}\label{eq:dvartwo}
	\V^2(X) =  T_1(X) +T_2(X) -2 \, T_3(X),
	\end{equation}
	where
	\begin{equation}\label{T1T2}
	\begin{aligned}
	T_1(X) &= \E (\|X-X'\|^2), \\
	T_2(X) &= (\E\|X-X'\|)^2, 
	\end{aligned}
	\end{equation}
	and 
	\begin{align}\label{T3}
	T_3(X) &= \E \big(\|X-X'\| \cdot\|X-X''\|\big).
	\end{align}
For an i.i.d. sample $\bX= (X_1,\ldots,X_n)$ drawn from $X$, the empirical version of $\V^2(X)$ was given in \citep{szekely2007} as
\begin{equation}\label{eq:sampledvar}
	\V_n^2(\bX) =  T_{1,n}(\bX) + T_{2,n}(\bX) -2 \, T_{3,n}(\bX),
	\end{equation}
	where 
	\begin{equation}\label{T1T2sample}
	\begin{aligned}
	T_{1,n}(\bX)&=\frac{1}{n^2}   \sum_{i=1}^n \sum_{j=1}^n \|X_i-X_j\|^2, \\
	T_{2,n}(\bX) &= \Big(\frac{1}{n^2} \sum_{i=1}^n \sum_{j=1}^n \|X_i-X_j\| \Big)^2, \\
	\end{aligned}
	\end{equation}
	and 
	\begin{align}\label{T3sample}
	T_{3,n}(\bX) &= \frac{1}{n^3}  \sum_{i=1}^n \sum_{j=1}^n \sum_{k=1}^n  \|X_i-X_j\| \cdot \|X_i-X_k\|.
	\end{align}
The version (\ref{eq:sampledvar}) is not unbiased; an unbiased estimator for $\V^2(X)$ was derived in \citep{Huo2015}, viz., 
	\begin{multline}
\widehat{\V}_n^2(\bX) = \frac{n}{n-3} \, T_{1,n}(\bX) + \frac{n^3}{(n-1) (n-2) (n-3)}  T_{2,n}(\bX) \\ %nonumber \\
-\frac{2 n^2}{(n-2) (n-3)}  \, T_{3,n}(\bX) \label{eq:Udcov}.
	\end{multline}
By \citep{Huo2015}, $\widehat{\V}_n^2(\bX)$ is a U-statistic of order four with kernel function
	\begin{align}
	\phantom{AAA} h(X_1,X_2,X_3,X_4) &= \frac{1}{4} \, \sum_{\substack{1 \leq i,j \leq 4 \\ i \neq j}} \|X_i - X_j \|^2 - \frac{1}{4} \, \sum_{i=1}^4 \Bigg(\sum_{\substack{j=1 \\ j \neq i}}^4 \|X_i - X_j \| \Bigg)^2 \nonumber \\
	& \qquad + \frac{1}{24} \Bigg(\sum_{\substack{1 \leq i,j \leq 4 \\ i \neq j}} \|X_i - X_j \| \Bigg)^2. \label{eq:kernel} 
	\end{align}
	
	In the sequel, we derive the asymptotic distribution of $\widehat{\V}_n^2(\bX)$ and  $\V_n^2(\bX)$; further, we do so under conditions weaker than known  previously. Hitherto, the asymptotic normality of $\widehat{\V}_n^2(\bX)$ was proved only under the assumption that the {\it fourth} moment of $X$ is finite; see \cite[Lemma 4.8 and Theorem 4.11]{Huang2017}. Here, we derive the asymptotic normality under the broader assumption that the second moment of $X$ is finite.
The following lemma provides an alternative representation for the kernel function $h(\cdot)$ in (\ref{eq:kernel}).
	\begin{lemma}\label{lem:kernel}
		The kernel function $h$ in (\ref{eq:kernel}) can be written as
		\begin{align*} 
		h(X_1,X_2,X_3,X_4) 
		&=\frac{1}{12} \, \sum_{\substack{1 \leq i,j \leq 4 \\ i \neq j}}^4  \|X_i - X_j \|^2 \\
		& \qquad - \frac{1}{12} \, \sum_{\substack{1 \leq i,j,k \leq 4 \\ i,j,k \ {\rm{ all \, different}}}} \|X_i - X_j \| \, \|X_i - X_k \| \\  
		& \qquad + \frac{1}{24} \sum_{\substack{1 \leq i,j,k,l \leq 4 \\ i,j,k,l \ {\rm{ all \, different}}}} \|X_i - X_j \| \,  \|X_k - X_l \|.
		\end{align*}
	\end{lemma}	
	
	In applying Lemma \ref{lem:kernel} to establish the asymptotic normality of $\V_n^2(\bX)$ under the assumption of finite second moments of $X$, let 
	$$
	h_1(x)=\E [h(x,X_2,X_3,X_4)] - \V^2 (X)
	$$ 
	be the linear part in the Hoeffding decomposition \cite[Section 11.4]{Vandervaart2000} of the kernel $h$ and let 
	\begin{equation} \label{eq:asv}
	\gamma = 16 \, \E [h_1^2(X)].
	\end{equation}
	
	We remark that an expansion of $h_1(X)$ in our setting is given in \cite[Eq. (B.6)]{Huang2017}.
	
	Denote by $\bX_{-k}$ the sample $\bX$ with the $k$-th observation deleted.  Assuming that $\E [h^2(X_1,X_2,X_3,X_4) ] < \infty$ (which is a consequence of $\E(\|X\|^2) < \infty$; cf. the proof of Theorem \ref{th:asymptotics}), it follows from Arvesen \cite[Theorem 9]{arvesen1969} that the jackknife estimator 
		\begin{equation}
			\widehat{\gamma}(\bX) = (n-1) \sum_{i=1}^n \Big( \widehat{\V}_{n-1}^2 (\bX_{-i}) - \frac{1}{n} \sum_{j=1}^n \widehat{\V}_{n-1}^2 (\bX_{-j}) \Big)^2 \label{eq:jkgamma}
		\end{equation}
	is a weakly consistent estimator of $\gamma$.

	\begin{theorem} \label{th:asymptotics}
		Suppose that $\E(\|X\|^2) < \infty$.  As $n \to \infty$, 
		\begin{align}
		\sqrt{n}\big(\widehat{\V}_n^2(\bX) - \V^2(X)\big) \cid N(0,\, \gamma) \label{Vnasymptoticdistn}
		\end{align}
		and
		\begin{align}
		\frac{\sqrt{n}\big(\widehat{\V}_n^2(\bX) - \V^2(X)\big)}{\sqrt{\widehat{\gamma}(\bX)}} \cid N(0,\, 1), \label{eq:Vnasymptoticgammahat}
		\end{align}
		and the same result holds for $\V_n^2(\bX)$. 
	\end{theorem}

	The asymptotic distribution of $\V_n(\bX)$, the empirical distance standard deviation, now follows from Theorem \ref{th:asymptotics} by the delta method. A weakly consistent estimator for the variance of the asymptotic distribution of $\sqrt{n} (\widehat{\V}_n(\bX) - \V(X))$ is obtained analogously from Arvesen \cite[Theorem 9]{arvesen1969} and is given by
	\begin{equation} \label{eq:jkxi}
			\widehat{\xi}(\bX) = (n-1) \sum_{i=1}^n \Big( \widehat{\V}_{n-1} (\bX_{-i}) - \frac{1}{n} \sum_{j=1}^n \widehat{\V}_{n-1} (\bX_{-j}) \Big)^2.
	\end{equation}
	
		\begin{corollary} \label{cor:asymptotics}
		Suppose that $\E(\|X\|^2) < \infty$ and $\V(X) > 0$. Then, 
		\begin{align*}
		\sqrt{n}(\widehat{\V}_n(\bX) - \V(X)) &\cid N\big(0, \gamma/4 \V^2(X)\big)
		\end{align*}
		and
		\begin{align*}
	\frac{\sqrt{n} (\widehat{\V}_n(\bX) - \V(X))}{\sqrt{\widehat{\xi}(\bX)}} &\cid N\big(0, 1\big),
		\end{align*}
		and the same result holds for $\V_n(\bX)$. 
	\end{corollary}

	We now consider the problem of estimating scale in a location-scale family of the form $X \stackrel{\tiny d}{=} \mu + \lambda Z$, with $\mu \in \R$ and  $\E |Z|^2 < \infty$, where $\stackrel{\tiny d}{=}$ denotes equality in distribution. In this location-scale setting, Corollary \ref{cor:asymptotics} enables the comparison of the efficiency of the distance standard deviation to other estimators of spread.  For any $\sqrt{n}$-consistent and asymptotically normal estimator $s_n(\bX)$, we define the asymptotic variance $\mbox{ASV}(s_n(\bX);F)$ at the distribution $F$ to be the variance of the limiting distribution of $\sqrt{n}(s_n(\bX) - s(X))$ as $n \to \infty$, where $s_n(\bX)$ is evaluated at an i.i.d.\ sequence drawn from $X \sim F$ and $s(X)$ denotes the corresponding population value of $s_n(\bX)$.
	 While two scale estimators (i.e., estimators satisfying property (C2) in Section 1) $s_n^{(1)}(\bX)$ and $s_n^{(2)}(\bX)$ may converge to different population values $s_1(X)$ and $s_2(X)$, respectively, $s_n^{(2)}(\bX)$ can be made consistent for $s_1(X)$ within the considered location-scale family by multiplying it with the factor $s_1(Z)/s_2(Z) = s_1(X)/s_2(X)$.  
	%Any estimators $s_n^{(1)}(\bX)$ and $s_n^{(2)}(\bX)$ that estimate %possibly different 
	%population values $s_1$ and $s_2$, respectively, at a given distribution $F$ and that satisfy the dilation property (C2) in Section 1 can be compared efficiency-wise by standardizing them by their respective population values.
	 Thus we define the {\it asymptotic relative efficiency} for scale estimators \cite[Eq. (2.1)]{Bickeldisp} $s_n^{(1)}(\bX)$ with respect to $s_n^{(2)}(\bX)$ at the population distribution 
	  $F$ as 
	\begin{equation} \label{eq:ARE}
	\mbox{ARE} \big(s_n^{(1)}(\bX),s_n^{(2)}(\bX);F\big) =  \frac{\mbox{ASV}(s_n^{(2)}(\bX);F)/(s_2(X))^2}{\mbox{ASV}(s_n^{(1)}(\bX);F)/(s_1(X))^2}.
	\end{equation}

		We consider as alternatives to $\V_n(\bX)$ the empirical standard deviation, 
	%$\widehat{\sigma}_n(\bX)$, which is defined as the positive square root of
		\begin{equation} \label{eq:empsd}
			\widehat{\sigma}_n(\bX) = \bigg[\frac{1}{n-1} \sum_{i=1}^n (X_i - \overline{X}_n)^2\bigg]^{1/2},
		\end{equation}
	where $\overline{X}_n$ denotes the sample mean of $\bX$, the empirical mean deviation
	  \begin{equation}
	  \widehat{d}_n(\bX) =  \frac{1}{n} \sum_{i=1}^n |X_i - m_n(\bX)|,
	  \end{equation}
	  	  where $m_n(\bX)$ denotes the sample median of $\bX$, and Gini's mean difference, 
		\begin{equation} \label{eq:empgini}
	  \widehat{\Delta}_n(\bX) = \frac{2}{n\,(n-1)} \sum_{1\leq i < j \leq n} |X_i - X_j|.
	 \end{equation}
We remark that$\big(\frac{n-1}{n} \big)^{1/2} \widehat{\sigma}_n(\bX)$ is the maximum likelihood estimator of scale in the location-scale family generated by the normal distribution $N(0,1)$.  Also, $\widehat{d}_n(\bX)$ is the analogous estimator of scale for the Laplace distribution $L(0,1)$.  

	\begin{table}[t!]
		\captionsetup{width=0.82\textwidth}
		\caption{\small The asymptotic relative efficiencies (\ref{eq:ARE}) with respect to the respective maximum likelihood estimators of the distance standard deviation $\V_n$, the standard deviation $\widehat{\sigma}_n$, the mean deviation $\widehat{d}_n$, and Gini's mean difference $\widehat{\Delta}_n$ at the Laplace distribution, the normal distribution, the normal scale mixture distribution $NM(3,0.1)$, and the $t_\nu$-distributions with $\nu = 3$ and $\nu = 5$. \label{tab:estimators.asymptotic} \label{TAB:ESTIMATORS.ASYMPTOTIC}}
		\renewcommand{\arraystretch}{1.1}
		\begin{center}
			\begin{tabular}{c|@{\quad } D{.}{.}{4}@{\qquad } D{.}{.}{4}@{\qquad } D{.}{.}{4}@{\qquad } D{.}{.}{4}@{\qquad } D{.}{.}{4}}  %-----------------------------------------------------------------------------------
				\hline
				Distribution, $F$  	{\phantom{\bigg|}}	&
				\multicolumn{1}{c}{\small $ARE(\V_n;F)$}	&
				\multicolumn{1}{c}{\small $ARE(\widehat{\sigma}_n;F)$}	&
				\multicolumn{1}{c}{\small $ARE(\widehat{d}_n;F)$}	&
				\multicolumn{1}{c}{\small $ARE(\widehat{\Delta}_n;F)$}	 \\
				\hline
				%-----------------------------------------------------------------------------------
				$L(0,1)$          & 0.952 & 0.8   & 1			& 0.964	\\
				$N(0,1)$          & 0.784 & 1     & 0.876 & 0.978	\\
				$N\!M(3,0.1)$	 & 0.887 & 0.398 & 0.757 & 0.641 \\
				$t_3$				 & 0.965 & 0		  &	0.681	& 0.524 \\
				$t_5$              & 0.992 & 0.4	  &	0.941	& 0.859 \\
				\hline 	
			\end{tabular}
		\end{center}
	\end{table}

Let $N\!M(\lambda,\epsilon)$ denote the normal scale mixture distribution that is defined as
\[
	N\!M(\lambda,\epsilon) 
	= (1 - \epsilon)N(0,1) + \epsilon N(0,\lambda^2), \qquad 0 \le \epsilon \le 1, \ \lambda \ge 1,
\]
and is also known as the \emph{contaminated normal distribution} \cite{tukey1960}.

	In Table \ref{tab:estimators.asymptotic}, we compare the asymptotic efficiencies of the distance standard deviation with the three alternative measures of spread at the Laplace distribution, normal distribution, the normal scale mixture distribution $N\!M(3,0.1)$, and the $t_\nu$-distributions with $\nu = 3$ and $\nu = 5$. 

The asymptotic relative efficiencies of these estimators with respect to the respective maximum likelihood estimator at each of the distributions are given in Table \ref{TAB:ESTIMATORS.ASYMPTOTIC}. 
Details on the calculations of the values in Table~\ref{tab:estimators.asymptotic} are given in Appendix \ref{appsec:table1} in the supplementary material. 
	
	While the distance standard deviation has moderate efficiency at normality, it turns out to be asymptotically very efficient in the case of heavier-tailed populations. For the normal scale mixture, the $t_3$- and the $t_5$- distributions, the distance standard deviation outperforms its three competitors. 

	\begin{table}[t]
		\captionsetup{width=0.8\textwidth}
		\caption{\small Simulated finite-sample values of the bias and the variance of the estimators $\V_n(\bX)$ and $\widehat{\V}_n(\bX)$ for $n = 5, 10, 50, 500$ compared to asymptotic values (last column); $10,000$ replications.
			\label{tab:estimators.finite}}
		\renewcommand{\arraystretch}{1.1}
		%\small
		\begin{center}
			\begin{tabular}{cc|@{\quad} D{.}{.}{4}@{\qquad } D{.}{.}{4}@{\qquad } D{.}{.}{4}@{\qquad } D{.}{.}{4}@{\qquad } D{.}{.}{4}}  %-----------------------------------------------------------------------------------
				Distribution     &  			& \multicolumn{5}{c}{Sample size} \\
				&	  &  \multicolumn{1}{l}{5} &  \multicolumn{1}{l}{10} & \multicolumn{1}{l}{50} & \multicolumn{1}{l}{500} & \multicolumn{1}{l}{$\infty$} \\
				\hline
				$L(0,1)$ & $\sqrt{n} (\E(\V_n)-\V)$ 					& 0.282 & 0.313 & 0.191 & 0.067 & 0 \\  
				& $\sqrt{n} (\E(\widehat{\V}_n)-\V)$ 				& -0.270 & -0.136 & -0.057 & -0.022 & 0\\ 
		& $n \Var(\V_n)$ 								& 0.953 & 0.834 & 0.668 & 0.605 & 0.613\\ 
	& $n \Var(\widehat{\V}_n)$ 			& 0.899 & 0.723 & 0.642 & 0.604 & 0.613\\[1.0ex] 
				$N(0,1)$ & $\sqrt{n} (\E(\V_n)-\V)$ 				&  0.067 & 0.085 & 0.049& 0.022 & 0\\   
				& $\sqrt{n} (\E(\widehat{\V}_n)-\V)$ 				& -0.197 & -0.082& -0.022 & 0.022 & 0\\
			& $n \Var(\V_n)$ 								& 0.299 & 0.280 & 0.256 & 0.253 & 0.256\\ 
	& $n \Var(\widehat{\V}_n)$ 			& 0.488 & 0.336 & 0.265 & 0.253 & 0.256\\[1.0ex]
				$N\!M(3,0.1)$ 	 & $\sqrt{n} (\E(\V_n)-\V)$ 	& 0.244 & 0.307 & 0.198& 0.067 & 0\\
				& $\sqrt{n} (\E(\widehat{\V}_n)-\V)$					& -0.246 & -0.114 & -0.049 &-0.022 & 0\\  
				& $n \Var(\V_n)$ 								& 0.905 & 0.807 & 0.514 & 0.440 & 0.426\\ 
			& $n \Var(\widehat{\V}_n)$ 			& 0.673 & 0.523 & 0.455 & 0.434 & 0.426\\[1.0ex] 
				$t_3$ 	 & $\sqrt{n} (\E(\V_n)-\V)$ 				& 0.409 & 0.471 & 0.368 & 0.157 & 0\\ 
				& $\sqrt{n} (\E(\widehat{\V}_n)-\V)$ 					& -0.304 & -0.158 & -0.049 &  0.000 & 0 \\
							& $n \Var(\V_n)$ 								& 4.231 & 2.027 & 1.123 & 0.746 & 0.680\\ 
			& $n \Var(\widehat{\V}_n)$ 			& 0.991 & 0.800 & 0.711 & 0.676 & 0.680\\[1.0ex] 
				$t_5$    & $\sqrt{n} (\E(\V_n)-\V)$ 			& 0.212 & 0.234 & 0.148 & 0.045 & 0\\ 
				& $\sqrt{n} (\E(\widehat{\V}_n)-\V)$ 					& -0.235 & -0.114 & -0.042 & -0.022 & 0\\ 
				& $n \Var(\V_n)$ 								& 0.772 & 0.638 & 0.472 & 0.427 & 0.424\\ 
			& $n \Var(\widehat{\V}_n)$ 			& 0.708 & 0.517 & 0.435 & 0.418 & 0.424\\[1.0ex]
				  \hline
			\end{tabular}
		\end{center}
	\end{table}

	In Table \ref{tab:estimators.finite}, we complement our asymptotic analysis with finite-sample simulations. For sample sizes $n = 5, 10, 50, 500$ and the same population distributions as above, the simulated biases and variances (based on $10,000$ replications) of the empirical versions of distance standard deviation $\V_n(\bX)$ and $\widehat\V_n(\bX)$ are given along with their respective asymptotic values.  The corresponding values for the competing estimators $\widehat{\sigma}_n(\bX)$, $\widehat{d}_n(\bX)$, and $\widehat{\Delta}_n(\bX)$ are provided by Gerstenberger and Vogel \citep[][Tables 7, 8]{gerstenberger2015}.  
	
	The values presented in Table \ref{tab:estimators.finite} indicate that $\widehat\V_n(\bX)$ is preferable to $\V_n(\bX)$ as an estimator of $\V(X)$.  We note that both estimators are biased and that $\widehat\V^2_n$ is a U-statistic whereas $\widehat\V_n$ is not; however, $\widehat\V_n$ shows considerably smaller bias than $\V_n$ at the heavier-tailed distributions. In light of the efficiency comparison with the other standard scale estimators, it emerges that heavy-tailed distributions represent the most promising area for applications of the distance standard deviation. 
	
	Although the definition of the empirical distance standard deviation $\widehat\V_n$ does not make apparent its superior performance under heavy tails, an intuitive explanation for its superiority in that context is obtained in Section \ref{sec:properties}, where the scale estimators are expressed in terms of the spacings between data points; it is seen there that the distance standard deviation provides comparably little weight to the extreme spacings at either end of the data range, and it is this property that leads to the superior performance of $\widehat\V_n$ with heavy-tailed data. 
	
	In concluding this section, we note that the main argument in the proof of Theorem \ref{th:asymptotics}, namely that $\E(\|X\|^2) < \infty$ implies	$\E[h^2(X_1,X_2,X_3,X_4)] < \infty$, leads to a proof of the central limit theorem (CLT) for the {\it squared distance covariance} given in \cite[Theorem 4.11]{Huang2017} under weaker conditions than known previously. The intrinsic idea in our proof of this CLT is that the square of the corresponding U-statistic,
		\begin{equation}
		\begin{aligned}
		&\widehat{\Omega}_n(\bX,\bY) = \frac{1}{n \, (n-3)}\Bigg[ \sum_{i,j=1}^n \|X_i-X_j\| \|Y_i-Y_j\| \\
		 &\qquad \qquad \qquad+  \frac{1}{(n-1) \, (n-2)} \sum_{i,j=1}^n \|X_i-X_j\| \cdot \sum_{i,j=1}^n \|Y_i-Y_j\| \\ 
		 &\qquad \qquad \qquad- \frac{2}{(n-2)} \sum_{i,j,k=1}^n \|X_i-X_j\| \|Y_i-Y_k\|\Bigg],\label{eq:ugendcov}
		\end{aligned}
		\end{equation}
is bounded above by $\widehat{\V}_n^2(\bX) \, \widehat{\V}_n^2(\bY)$, where $\bY = (Y_1,\ldots,Y_n)$ denotes an i.i.d. sample drawn from some random variable $Y \in \R^q$. The complete statement of this limit theorem is given in Section \ref{app:clt} of the supplementary material.

	\section{Applications of the Distance Variance}{\label{sec:app}}

	We consider two applications in detail: hypothesis testing and quality control charts.  For other interesting applications of the distance standard deviation, we refer to Fiedler \cite{fiedler2016a} who defined the distance variogram and gave a natural generalization of the usual variogram for $\alpha$-stable distributions.
	
	\subsection{Hypothesis Testing}\label{subsec:hyptesting}

	For ease of exposition, we focus on two-sample hypothesis testing in the univariate case. One-sample tests and results for the multivariate setting can be derived analogously.

	Let $\bX_n = (X_1,\ldots,X_n)$ and $\bY\!_m = (Y_1,\ldots,Y_m)$ be two i.i.d., mutually independent random samples drawn from random variables $X,Y \in \R$ with finite second moments. We wish to test the null hypothesis $H_0: \V(X) = \V(Y)$.  
	For this purpose, we propose the test statistic 
	\begin{equation} \label{eq:dsdtest}
	\widehat{T}_\V = \sqrt{\frac{n \, m}{n + m}} \, \frac{\widehat\V_n(\bX_n) - \widehat\V_m(\bY\!_m)}{ \sqrt{\widehat{\xi}_p(\bX_n,\bY\!_m)}}, 
	\end{equation}
	where $\widehat{\xi}_p(\bX_n,\bY\!_m)$ is a pooled estimator of the form
	$$
	\widehat{\xi}_p(\bX_n,\bY\!_m) = \frac{n \, \widehat{\xi}(\bX_n) + m \, \widehat{\xi}(\bY\!_m)}{n +m}
	$$
and $\widehat{\xi}(\cdot)$ is defined in equation (\ref{eq:jkxi}). By Theorem \ref{th:testlim}, $\widehat{T}_\V$ and $|\widehat{T}_\V|$ can be directly applied to test $H_0$ against one-sided and two-sided alternatives, respectively. In stating this theorem, we denote by $\xi_X$ and $\xi_Y$ the asymptotic variances of the distributions of $\sqrt{n} (\widehat\V_n(\bX_n)-\V(X))$ and $\sqrt{m} (\widehat\V_m(\bY_m)-\V(Y))$, respectively (see Corollary \ref{cor:asymptotics} for details). %Moreover, for two random variables $Z,U$, we write $Z \stackrel{d}{=} U$, if $Z$ and $U$ follow the same distribution. 

\begin{theorem} {\label{th:testlim}}
	 Let $\E|X|^2 < \infty$ and $\E|Y|^2 < \infty$. Then, for $n,m \to \infty$, such that $n/m \to r >0$, it holds % $\widehat{T}_\V$ shows the following limit behavior:
	 \begin{enumerate}
	 	\item[(i)] If $\V(X) = \V(Y)$ then $\widehat{T}_\V \cid N\big(0,(\xi_X+r \xi_Y)/(r \xi_X + \xi_Y)\big)$.
	 	In particular, if additionally $X + \mu \stackrel{d}{=} Y$ where $\mu \in \R$, or $n/m \to 1$, then
	 	$\widehat{T}_\V \cid N(0,1)$. 
	 \item[(ii)] If $\V(X) < \V(Y)$ then $\widehat{T}_\V \cip -\infty$.
	  \item[(iii)] If $\V(X) > \V(Y)$ then $\widehat{T}_\V \cip \infty$.
		\end{enumerate}
	\end{theorem}
	
	A weakly consistent estimator of the asymptotic variance in Theorem \ref{th:testlim}(i) is 
	$$
	\big[\widehat{\xi}(\bX_n)+(n/m)\widehat{\xi}(\bY_m)\big]{\big/}\big[(n/m)\widehat{\xi}(\bX_n)+\widehat{\xi}(\bY_m)\big],
	$$
	and this estimator can be used to construct a Studentized statistic for testing $H_0: \V(X) = \V(Y)$.  For cases in which $n/m \to 1$ or if $X$ and $Y$ belong to a common location-scale family, the resulting asymptotic variance equals 1, identically, and hence no estimation is needed.  
	
	When the distributions of $X$ and $Y$ belong to the same location-scale family, i.e., $\lambda X + \mu \stackrel{d}{=} Y$  with $\mu \in \R$, the null hypothesis can be expressed as $H_0:  \lambda = 1$. Within this location-scale setting, we can compare the distance standard deviation based two-sample test to analogously constructed tests based on the standard deviation and Gini's mean difference. For scale measures $s_1$ and $s_2$, we denote by ${s_n^{(1)}}(\cdot)$ and ${s_n^{(2)}}(\cdot)$ the respective empirical versions of these measures. Moreover, we assume for i.i.d. samples $\bX_n$ drawn from random variables $X \sim F$ and for $j = 1,2$ that:
	\begin{itemize}
		\item[(A1)] 	If $y_i = b x_i + a$ for $a,b \in \R, i \in \{1,\ldots,n\}$, $n \in \N$, then ${s_n^{(j)}}(\by\!_n) = |b| \, {s_n^{(j)}}(\bx_n)$, where $\bx_n = (x_1,\ldots,x_n)$ and $\by_n = (y_1,\ldots,y_n)$.
		\item[(A2)] For $n \to \infty$, $\sqrt{n} \, ({s_n^{(j)}}(\bX_n)- s_j (X) ) \stackrel{d}{\longrightarrow} \mathcal{N}(0,\xi_j)$ with $\xi_j > 0$ and $\widehat{\xi}_j(\bX_n)$ is a consistent estimator for $\xi_j$.
		\item[(A3)]  The estimator of the asymptotic standard deviation $\sqrt{\widehat{\xi}_j}$ satisfies (A1) (with $ {s_n^{(j)}}$ replaced by  $\sqrt{\widehat{\xi}_j}$).
	\end{itemize}
	Test statistics analogous to (\ref{eq:dsdtest}) can then be constructed as
	\begin{equation} \label{eq:test statistic}
	\widehat{T}_j(\bX_n, \bY\!_m) = \sqrt{\frac{n \,m}{n+m}} \frac{{s_n^{(j)}}(\bX_n) - s_n^{(j)}(\bY\!_m)}{\sqrt{\widehat{\xi}_{p,j}(\bX_n,\bY\!_m)}},
	\end{equation}
	where
	\[
	\widehat{\xi}_{p,j}(\bX_n,\bY\!_m) = \frac{n}{n+m} \widehat\xi_{j}(\bX_n) + \frac{m}{n+m} \widehat\xi_{j}(\bY\!_m). 
	\]
	Theorem \ref{th:ARE} provides a comparison of the efficiency of two scale tests of the form (\ref{eq:test statistic}) under local alternatives. 
		Let $\lambda_{n,m}$ be an array of real numbers satisfying 
	$$
	\sqrt{\frac{n \,m}{n+m}} (\lambda_{n,m} - 1)  \to \Lambda
	$$
	for some $\Lambda \in \R$ as $n,m\to\infty$ such that $n/m \to r > 0$.
	For $m,n \in \N$, denote by $\bX_n = (X_1,\ldots,X_n)$, $\bZ_m = (Z_1,\ldots,Z_m)$ two mutually independent, i.i.d. samples  drawn from random variables $X$ and $Z$ following the same distribution $F$. Moreover, for $k \leq m$, define
	\begin{equation} \label{eq:localalt}
	\bY^{(n,m)}_k = (\lambda_{n,m} Z_1 + \mu, \ldots, \lambda_{n,m} Z_k + \mu),
	\end{equation}
	with $\mu \in \R$. In the following theorem, $[t]$ will denote the integer part of $t \in \R$ and 		$\Phi^{-1}(\cdot)$ is the inverse of $\Phi$, the standard normal distribution function.

	\begin{theorem} \label{th:ARE}
	 Let  $\rho = \xi_2 s_1^2(X)/ (\xi_1 s_2^2(X))$ denote the asymptotic relative efficiency of ${s_n^{(1)}}$ with respect to ${s_n^{(2)}}$ at $F$, cf.~(\ref{eq:ARE}), where we assume without loss of generality that $\rho \leq 1$. 
		Then, under Assumptions (A1), (A2), and (A3), 
		$\widehat{T}_1(\bX_n, \bY_m^{(n,m)})$ and $\widehat{T}_2(\bX_{[\rho n]}, \bY_{[\rho m]}^{(n,m)})$
		both converge in distribution to $N(- s_1(X)\Lambda/\sqrt{\xi_1}, 1)$ as $n, m \to \infty$, such that $n/m \to r > 0$. 
		
		Consequently, for the ratio of the power of two-sided tests with asymptotic size $\alpha$,
		\[
		%\mbox{\Large $\displaystyle\lim_{\substack{n,m \to \infty \\ n/m \to r}}$}
		\lim_{\substack{n,m \to \infty \\ n/m \to r}}
		\frac{P\left( \big|\widehat{T}_1(\bX_n, \bY_m^{(n,m)})\big| > \Phi^{-1}(1- \alpha/2)  \right)}{P\left( \big|\widehat{T}_2(\bX_{[\rho n]}, \bY_{[\rho m]}^{(n,m)})\big| > \Phi^{-1}(1- \alpha/2) \right)} = 1.
		\]
	\end{theorem}

	Two benefits of Theorem \ref{th:ARE} are that it enables explicit calculation of the asymptotic power for alternatives of the form (\ref{eq:localalt}), and it establishes a direct link between the asymptotic relative efficiencies of the scale estimators studied in Section \ref{sec:estimator} and the efficiencies of corresponding two-sample tests.
	
Let $\widehat{T}_\sigma$ and $\widehat{T}_\Delta$ denote test statistics according to (\ref{eq:test statistic}) based on 
$\widehat\sigma_n$ and 
$\widehat\Delta_n$, respectively, where the asymptotic variance of each estimator is estimated by the jackknife method, cf.~(\ref{eq:jkxi}). 
\begin{table}[t!]

	\captionsetup{width=\textwidth}
	\caption{\small  \emph{Test size.}  Empirical rejection frequencies (\%) under the null hypothesis $\lambda=1$ of asymptotic two-sample scale tests (based on the distance standard deviation $\widehat\V_n$, the standard deviation $\widehat\sigma_n$, Gini's mean difference $\widehat\Delta_n$, and the $F$-test) at the 5\% significance level. Results are based on $10,000$ replications.
		\label{tab:tests.size}}
	\renewcommand{\arraystretch}{1.05}
	%\small
	\begin{center}
		\begin{tabular}{cc|ccccccc} 
			%----------------------------------------------------------------------------------
			& $n$ 				& 15 & 50 & 120 	& 250 & 600 & 1,000 & $\infty$ \\
	& $m$ & 15  & 50 	& 40 		& 250 		& 200 		& 1,000 			& $\infty$ \\[1.1ex]
			Distribution & Test & \multicolumn{5}{c}{Rejection frequencies (\%)} \\[1.2ex]
			\hline {\phantom{\Big[}}
			$L(0,1)$ 	& $\widehat\V_n$  		& 4.4   & 4.6  & 5.0 	& 4.8 & 5.2 	& 4.7  & 5.0 \\
			& $\widehat\sigma_n$ 	& 3.3   & 4.1  & 4.4 	& 4.6 & 4.9 	& 5.0  & 5.0\\
			& $\widehat\Delta_n$ 	& 5.9   & 5.1  & 4.9 	& 4.8 & 5.3 	& 4.9  & 5.0\\
			& $F$-test 				& 17.5  & 20.0 & 19.5 & 21.6 & 20.8 & 20.6  & \\[1.5ex] 
			$N(0,1)$ 	& $\widehat\V_n$ 			& 4.5  	& 5.1 & 5.1 	& 5.3 & 5.1 & 5.2  & 5.0\\
			& $\widehat\sigma_n$ 	& 4.4 	& 5.0 & 5.0 	& 5.4 & 5.1 & 5.1  & 5.0\\
			& $\widehat\Delta_n$	& 5.8  	& 5.3 & 5.3 	& 5.4 & 5.3 & 5.0  & 5.0\\
			& $F$-test 				& 5.4  	& 5.3 & 4.9 	& 5.3 & 5.2 & 5.1  & 5.0\\[1.5ex] 
			$NM(3,0.1)$ & $\widehat\V_n$ 		& 4.0 	& 4.2 & 4.3 	& 4.9 & 4.8 & 5.0  & 5.0\\
			& $\widehat\sigma_n$ 	& 2.6 	& 3.3	& 4.3 	& 4.4 & 4.8 & 4.9  & 5.0\\
			& $\widehat\Delta_n$ 	& 4.7 	& 5.2 & 4.8 	& 4.8 & 5.0 & 5.0  & 5.0\\
			& $F$-test 				& 21.0 	& 27.3 & 27.9	& 30.6& 31.3& 31.1 & \\[1.5ex] 
			$t_3$ 		& $\widehat\V_n$ 			& 4.3  & 4.5 & 4.9 		& 4.6 & 4.8 & 5.2  & 5.0\\
			& $\widehat\sigma_n$ 	& 3.0  & 2.7 & 4.1 		& 3.3 & 4.3 & 3.5  &    \\
			& $\widehat\Delta_n$ 	& 5.4  & 4.6 & 5.0 		& 4.7 & 4.7 & 4.7  & 5.0\\
			& $F$-test 				& 25.8  & 35.8 & 38.0 & 49.9& 51.2& 59.3 &    \\[1.5ex] 
			$t_5$ 		& $\widehat\V_n$ 			& 4.4  & 4.8 & 4.7 		& 4.7 & 4.8 & 5.2  & 5.0\\
			& $\widehat\sigma_n$ 	& 3.6  & 4.0 & 4.4 		& 4.4 & 4.9 & 4.4  & 5.0\\
			& $\widehat\Delta_n$ 	& 5.6  & 5.2 & 5.0 		& 5.0 & 5.0 & 5.1  & 5.0\\
			& $F$-test 				& 14.5 & 19.4 & 18.2 	& 24.2& 24.5& 26.7 &    \\[1.5ex]  
			\hline
		\end{tabular}
	\end{center}
\end{table}	

\begin{table}[t!]
	\captionsetup{width=\textwidth}
	\caption{\small \emph{Test power.} Empirical rejection frequencies (\%) under the alternative $\lambda_{n,m} = 1 + 3 \sqrt{(n+m)/n/m}$ of asymptotic two-sample scale tests (tests based on the distance standard deviation $\widehat\V_n$, the standard deviation $\widehat\sigma_n$, Gini's mean difference $\widehat\Delta_n$, and the $F$-test) at the 5\% significance level. Results are based on $10,000$ replications.		
		\label{tab:tests.power}}
	\renewcommand{\arraystretch}{1.1}
	%\small
	\begin{center}
		\begin{tabular}{cc|ccccccc} 
			%----------------------------------------------------------------------------------
			& $n$ 				& 15 & 50 & 120 	& 250 & 600 & 1,000 & $\infty$ \\
			& $m$ & 15  & 50 	& 40 		& 250 		& 200 		& 1,000 			& $\infty$ \\[1.1ex]
			Distribution & Test & \multicolumn{5}{c}{Rejection frequencies (\%)} \\[1.1ex]
			\hline     {\phantom{\Big[}}
			$L(0,1)$ 	& $\widehat\V_n$  		& 33.0  & 57.1 & 63.8 & 72.6 & 75.0 & 78.8 & 83.3 \\
			& $\widehat\sigma_n$ 	& 27.1 & 51.5 & 66.1 & 65.8 & 71.5 & 71.7 & 76.5 \\
			& $\widehat\Delta_n$ 	& 44.6 & 62.2 & 68.7 & 73.6 & 76.6 & 79.3 & 83.8 \\[1.5ex] 
			$N(0,1)$ 	& $\widehat\V_n$ 			& 48.0 & 78.0 & 83.2 & 90.9 & 91.8 & 94.0 & 96.4 \\
			& $\widehat\sigma_n$ 	& 56.9 & 87.3 & 92.5 & 96.0 & 97.0 & 97.9 & 98.9 \\
			& $\widehat\Delta_n$	& 68.7 & 88.1 & 91.4 & 95.8 & 96.3 & 97.6 & 98.7 \\
			& $F$-test 				& 76.1 & 90.3 & 92.0 & 96.2 & 96.6 & 97.9 &       \\[1.5ex] 
			$NM(3,0.1)$ & $\widehat\V_n$ 		& 41.5 & 68.2 & 73.1 & 82.7 & 83.9 & 87.2 & 91.1 \\
			& $\widehat\sigma_n$ 	& 28.2 & 40.2 & 55.0 & 51.0 & 56.2 & 54.7 & 60.0 \\
			& $\widehat\Delta_n$ 	& 46.3 & 60.2 & 66.0 & 70.8 & 73.0 & 74.8 & 80.3 \\[1.5ex] 
			$t_3$ 		& $\widehat\V_n$ 			& 36.3 & 58.3 & 64.4 & 74.2 & 75.2 & 79.0 & 83.8 \\
			& $\widehat\sigma_n$ 	& 23.1 & 31.9 & 45.1 & 31.3 & 35.8 & 24.8 &  			\\
			& $\widehat\Delta_n$ 	& 39.9 & 49.9 & 55.6 & 56.4 & 58.1 & 58.1 & 58.4 \\[1.5ex] 
			$t_5$ 		& $\widehat\V_n$ 			& 41.6 & 69.0 & 73.9 & 82.8 & 85.2 & 88.4 & 91.6 \\
			& $\widehat\sigma_n$ 	& 36.4 & 56.5 & 68.7 & 64.1 & 68.2 & 64.7 & 56.4 \\
			& $\widehat\Delta_n$ 	& 53.6 & 69.2 & 74.1 & 78.6 & 81.4 & 83.2 & 87.5 \\	[1.5ex]  		
			\hline
		\end{tabular}
	\end{center}
\end{table}
%
%% description of tables 	
Tables \ref{tab:tests.size} and \ref{tab:tests.power} contain rejection frequencies (based on $10,000$ replications) at the $5\%$ level for two-sided asymptotic tests based on $|\widehat{T}_\V|$, $|\widehat{T}_\sigma|$ and $|\widehat{T}_\Delta|$. % using the asymptotic normal approximation and the jackknife method for variance estimation.
The $F$-test is also included for the sake of completeness and to serve as a benchmark in the normal case. The sensitivity of the $F$-test with respect to the assumption of normality is well known and is confirmed by the tables. 
We consider the Laplace distribution, normal distribution, normal scale mixture distribution $NM(3,0.1)$, and the $t_\nu$-distributions with $\nu = 3$ and $\nu = 5$.
The sample sizes $n, m$ range from $n + m = 30$ to $n + m = 2,000$. Table \ref{tab:tests.size} (test size) contains results for the null hypothesis $\lambda = 1$ and Table \ref{tab:tests.power} (test power) gives results for the sample-size-dependent alternative with
\begin{equation} \nonumber	%\label{lambdanm}
\lambda_{n, m} = 1 + 3\sqrt{\frac{n+m}{n \, m}}. 
\end{equation}
Theorem \ref{th:ARE}
 yields large-sample approximations for the power of the tests, which are provided in the last column of Table 3 and Table 4. The asymptotic power for the distance standard deviation test is $P(|N(3 \V /\sqrt{\xi}, 1)| > \Phi^{-1}(0.975))$, and similar expressions hold for $\widehat \sigma_n$ and $\widehat \Delta_n$. Note that $\widehat\sigma_n$ does not satisfy the conditions of Theorem \ref{th:ARE} at the $t_3$-distribution. 
%
%
%

%% interpretation of results of tables
%
In Table \ref{tab:tests.size}, we observe that the tests $|\widehat{T}_\V|$, $|\widehat{T}_\sigma|$, and $|\widehat{T}_\Delta|$ control the nominal level of 5\% well for all distributions under consideration. The actual rejection frequencies for the distance standard deviation test $|\widehat{T}_\V|$ range between 4.0 and 5.3. 
The $F$-test grossly exceeds the nominal level for non-normal distributions and is therefore omitted from the power considerations in Table \ref{tab:tests.power} except for the normal case. 

In Table \ref{tab:tests.power} we see that except for the small-sample case $(n,m)=(15,15)$, the distance standard deviation test $|\widehat{T}_\V|$ performs best at the heavier-tailed distributions $NM(3,0.1)$, $t_3$, and $t_5$. At the Laplace $L(0,1)$, $|\widehat{T}_\Delta|$ performs best and outperforms $|\widehat{T}_\V|$ and $|\widehat{T}_\sigma|$ for small sizes; for large sample sizes $|\widehat{T}_\V|$ and $|\widehat{T}_\Delta|$ perform almost equally. At the normal distribution $N(0,1)$, $|\widehat{T}_\sigma|$ and $|\widehat{T}_\Delta|$ dominate $|\widehat{T}_\V|$.

For small sample sizes, a better performance of the considered two-sample scale tests may be achieved by using a permutation-based approach for obtaining critical values, which we investigate in Section \ref{app:permtest} in the supplementary material.  We note that the permutation test requires both distributions to share a common location, which is a more restrictive assumption than is needed for the asymptotic test.

Finally, we remark that, in the univariate case, the distance variance and hence the distance standard deviation can be computed rapidly. For the asymptotic derivations in Section \ref{sec:estimator}, we used a fourth-order U-statistic representation of $\widehat\V_n^2$, which may suggest the opposite; however, Huo and Sz\'ekely \cite{Huo2015} devised an $O(n \log n)$ algorithm, which shows that the distance standard deviation has the same computational complexity as Gini's mean difference. All calculations for this article have been carried out using computationally efficient implementations of the distance standard deviation from the R package \texttt{dcortools} \cite{dcortools}, available on \url{https://github.com/edelmand21/dcortools}. An alternative $O(n \log n)$ implementation for the distance standard deviation is provided in the R package \texttt{energy} \cite{energy}, available on the Comprehensive R Archive Network (CRAN). 
	\paragraph{Data Examples}

\begin{figure}[t!]
		\centering
		\captionsetup{width=0.8\textwidth}
		\includegraphics[width=1.0\textwidth]{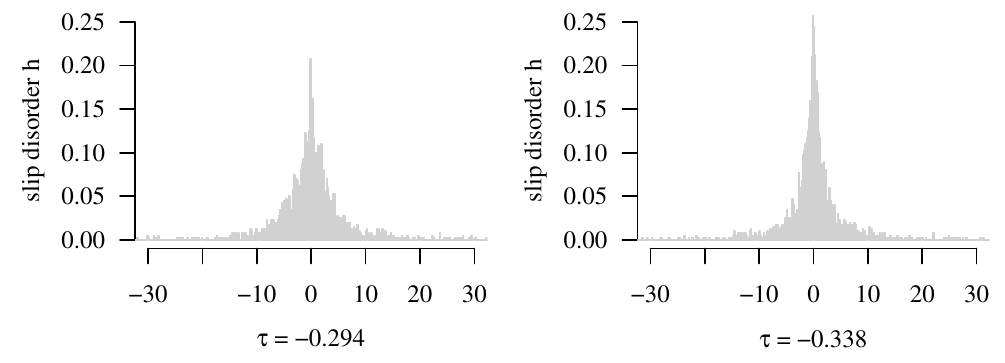}
		\caption{\small Histograms of slip disorder for two values of the thermal driving parameter $\tau$; $n = 2,000$ observations each.}
		\label{fig:data.ex.1}
	\end{figure}

We demonstrate the use of the two-sample distance variance test with an application to a data example. The data set stems from a physical model, studied by \citet{perez-reche:2016}, for describing a crystal undergoing a structural phase transition between austensite and martensite phases when subject to several cooling-heating cycles. In this model a quantity called the \emph{slip disorder}, and denoted by $h$, is of particular interest. The slip disorder depends on a parameter $\tau$ that represents thermal driving and is referred to as the \emph{temperature} within the model. The distribution of $h$ for two values of $\tau$ is depicted in Figure~\ref{fig:data.ex.1}.

An important issue is whether, and how, the distribution of the slip disorder $h$ is affected by $\tau$.  
As a consequence of the theoretical properties of the model, the distribution of $h$ is symmetric around zero; this symmetry is also suggested by Figure~\ref{fig:data.ex.1}. Hence, the main feature by which $h$ may differ with respect to $\tau$ is in its scale. The distribution of observed values of $h$ is very heavy-tailed and clearly non-normal, with excess kurtoses exceeding 10 (cf.~\cite[Figures 10 and 11]{perez-reche:2016}).  Consequently, an $F$-test is inappropriate here. 

	Although arbitrarily large data sets can be obtained by letting the model run sufficiently long, the simulations are computationally costly. Hence, fast detection of a statistically significant difference is desirable. For the data depicted in Figure~\ref{fig:data.ex.1} (with parameter values $\tau = -0.294$ and $\tau = -0.338$, and sample sizes $2,000$ each), the two-sided asymptotic test based on the statistic $|\widehat{T}_\V|$ yields a $p$-value of $0.0013$. Corresponding tests based on the standard deviation and Gini's mean difference give $p$-values of 0.6940 and 0.0327, respectively. 
The large $p$-value of the test based on the classical standard deviation is consistent with the tendency of the test to under-reject the null hypothesis in the presence of heavy tails; see Table \ref{tab:tests.power}.
	Moreover, as the standard deviation test remains persistently non-significant even for samples of size $10,000$, we find again that the classical standard deviation is an inappropriate measure of spread for heavy-tailed distributions. 	

	\subsection{Multivariate Statistical Quality Control}\label{subsec:MSQC}
	
	In statistical quality control \cite{Chandra2001,montgomery2009}, the objective is to monitor quality characteristics in production processes using statistical methods.  Important tools for process monitoring are the Shewhart control charts that are used to survey whether quality characteristics of the production process are under control. 
	
	A common tool used to monitor the dispersion of multivariate processes is the $|S|$-chart \cite{aparisi1999}, which applies the {\it generalized variance}, i.e., the determinant of the covariance matrix of the process. In the following, we investigate the potential of a control chart based on $\widehat{\V}_n$ as an alternative tool for monitoring the dispersion of multivariate processes. The $\V$-chart, as we will call the corresponding chart, can also be applied in high-dimensional settings, where the dimension of the production process exceeds the number of samples per subgroup (as for example in molecular %high-throughput
	 data).

	 To compare the performance of the $|S|$-chart with the $\V$-chart, we conduct a simulation study. For simplicity, we assume throughout the simulation study that the process under consideration is bivariate and that the two components of the production process are independent. The covariance matrix of the process in control will always be given by
	$$
	\Sigma_0 = \begin{pmatrix}
	1 &0 \\ 0 & 1
	\end{pmatrix}.
	$$
	For the covariance matrix of the process out of control, we will consider the matrices
	$$
	\Sigma_1 = \begin{pmatrix}
	\delta^2 &0 \\ 0 & 1
	\end{pmatrix}, \quad \quad \Sigma_2 = \begin{pmatrix}
	\delta &0 \\ 0 & \delta
	\end{pmatrix},
	$$		
	where $\delta=1,1.5,2,2.5,3,3.5,4$ and the distribution of the components will follow either a normal, Laplace, $t_3$- or $t_5$-distribution. 
	A Shewhart control chart consists of an upper control limit (UCL) and a lower control limit (LCL) for the characteristic under consideration. 
	This characteristic (e.g. the generalized variance) is then computed for consecutive subgroups of a fixed size $k$. When the characteristic lies below the LCL or exceeds the UCL, this represents an out-of-control signal of the process; the corresponding subgroup can then be investigated and, if necessary, further action can be taken. Formally, this corresponds to consecutive testing of the null hypothesis that the characteristic in a subgroup under consideration equals the characteristic in-control. 
	
	To compare the $|S|$- and $\V$-charts, we pursue a bootstrap approach. First, we simulate $10,000$ i.i.d. samples of the process in control, i.e. using the covariance matrix $\Sigma_0$. In application, this is typically given by a phase-I sample of the process, for which it is known that the production process was in control. From these $10,000$ samples, we now take $B=100,000$ bootstrap samples of size $k=25$, where $k$ coincides with the fixed subgroup size. For each of the bootstrap samples, the generalized variance and the distance standard deviation are evaluated. The UCL for the $|S|$-chart is then given by the $99.75 \%$-quantile of the generalized variances of the $B$ bootstrap samples, the corresponding LCL is given by the $0.25 \%$-quantile. The respective UCL and LCL for the $\V$-chart are calculated analogously. Using the respective alternative distribution (i.e. using $\Sigma_1$ or $\Sigma_2$), we now generate $250,000$ i.i.d. samples which are partitioned into $10,000$ subgroups of size $25$. For each subgroup, we evaluate if the generalized variance (or distance standard deviation respectively) exceeds the bounds given by the UCL or LCL of the respective chart. The empirical power for the $|S|$-chart and the $\V$-chart are then calculated by the fraction of subgroups for which these bounds were exceeded. The procedure is replicated $N=100$ times and the empirical power is averaged over these runs. The two methods are then compared using the average run length (ARL), which is the average number of subgroups one needs to test until an out-of-control signal is obtained. In the case of i.i.d. samples, the ARL is given by the reciprocal value of the empirical power.

 	\begin{table}[t]
 	\captionsetup{width=0.75\textwidth}
 	\caption{\small Average run lengths for the $\V$-chart and $|S|$-chart, respectively, for several alternative distributions representing the out-of-control state.}
 	\renewcommand{\arraystretch}{1.1}
 	%\small
 	\begin{center}
 		\begin{tabular}{cc|ccccccc} 
 			%-----------------------------------------------------------------------------------
 		Distribution  & Chart					& \multicolumn{7}{c}{Effect size, $\delta$} \\
 			&           		&  1 & 1.5 & 2  & 2.5  & 3  & 3.5 & 4  \\
 			\hline					
 		normal  ($\Sigma_1$) 	& $\V$ & 198.77	& 2.72		& 1.19 &  1.01 & 1.00 & 1.00 & 1.00     \\ 
 			& $|S|$	&	198.77			& 3.77 & 1.30 &	1.05	 & 1.01 & 1.00 & 1.00    \\ 
normal  ($\Sigma_2$) 	& $\V$ & 202.92		& 6.11 &  1.70 & 1.15 & 1.03 & 1.01 & 1.00     \\ 
& $|S|$					&	200.48			& 3.79 & 1.31 &	1.05	 & 1.01 & 1.00 & 1.00   \\ [1.0ex] 
 			Laplace  ($\Sigma_1$) 	& $\V$ 	&	196.23	& 5.85 &  1.53 & 1.10 & 1.02 & 1.00 & 1.00     \\ 
 		& $|S|$	&	190.30			& 13.16 & 3.18 &	1.69 & 1.27 & 1.12 & 1.05    \\ 
 		Laplace  ($\Sigma_2$) 	& $\V$ & 199.80		& 14.07 &  3.38 & 1.75 & 1.30 & 1.13 & 1.06     \\ 
 		& $|S|$					& 193.27 & 13.21 & 3.24	& 1.69 & 1.27 & 1.11 & 1.05   \\ [1.0ex] 
 		 			$t_5$  ($\Sigma_1$) & $\V$ 	& 205.38 & 4.70 & 1.32	& 1.03 & 1.00 & 1.00 & 1.00      \\ 
 		& $|S|$	&	184.20			& 26.99 & 5.99 &	2.55				 & 1.59 & 1.27 & 1.15    \\ 
 		$t_5$   ($\Sigma_2$) 	& $\V$  	& 199.23		& 11.68 & 2.69 & 1.47 & 1.16 & 1.05 & 1.02     \\ 
 		& $|S|$					& 185.77 & 29.63 &	6.19 & 2.62 & 1.65 & 1.28 & 1.14   \\ [1.0ex] 
 			$t_3$  ($\Sigma_1$) 	& $\V$  	& 199.80 & 12.78 & 2.20	& 1.22 & 1.04 & 1.01 & 1.00     \\ 
 		& $|S|$	&	169.00			& 94.04 & 40.40 &	18.41				 & 9.17 & 5.87 & 4.51    \\ 
 		$t_3$   ($\Sigma_2$) 	& $\V$  	& 202.92		& 34.12 & 7.01 & 2.99 & 1.80 & 1.37 & 1.17      \\ 
 		& $|S|$					& 176.37 & 84.69& 33.07 & 18.00 & 10.45 & 6.14 & 4.41   \\
 			\hline
 		\end{tabular}
 	\end{center}
 	\label{tab:runlength}
 \end{table}

	Table \ref{tab:runlength} lists the ARLs of the $|S|$-chart and the $\V$-chart for each of the covariance matrices $\Sigma_1$ and $\Sigma_2$, where $\delta=1,1.5,2,2.5,3,3.5,4$. The $\V$-chart shows substantial advantages compared to the $|S|$-chart for heavy-tailed distributions, such as the $t_3$- and $t_5$-distributions. Moreover, while the $|S|$-chart shows comparable performances for the different dispersion settings given by $\Sigma_1$ and $\Sigma_2$ (which is not surprising since they feature the same generalized variance), the distance standard deviation seems to be more powerful against large changes in one component compared to moderate changes in both components. For a definitive statement about the potential of the distance standard deviation for multivariate statistical quality control, more detailed comparisons with the generalized variance are required that go beyond the scope of this paper. Yet, our results indicate that the $\V$-chart is a promising alternative to the generalized variance for multivariate statistical quality control in the presence of heavy-tailed distributions.

		\section{Inequalities between the distance variance, the variance, and Gini's mean difference}\label{sec:ineq}

	In the following we will study inequalities between the summands appearing in (\ref{eq:dvartwo}) and (\ref{eq:sampledvar}). In the one-dimensional case, these inequalities will lead to crucial results concerning the relationships between $\V(X)$, $\Delta(X)$, and $\sigma(X)$.  
	
		\begin{lemma}{\label{lem:Tineq}}
		Let $T_{1,n}(\bX)$, $T_{2,n}(\bX)$, $T_{3,n}(\bX)$ be defined as in (\ref{T1T2sample}) and (\ref{T3sample}). Then there hold the algebraic inequalities, 
		\begin{equation}\label{Tineqs1and2}
		T_{2,n}(\bX) \leq  T_{3,n}(\bX) \leq T_{1,n}(\bX), \quad \quad \quad T_{1,n}(\bX) \leq 2 \, T_{3,n}(\bX).
		\end{equation}
		Further, if $X \in \R^p$ is a random vector such that $\E\|X\|^2 < \infty$, and if $T_1(X)$, $T_2(X)$, $T_3(X)$ are defined as in (\ref{T1T2}) and (\ref{T3}) then,
		\begin{equation}\label{Tineqs3and4}
		T_2(X) \leq  T_3(X) \leq T_1(X), \quad \quad \quad  T_1(X) \leq 2 \, T_3(X).
		\end{equation}
	\end{lemma}
	
%	\medskip
	
	Using the inequalities in Lemma \ref{lem:Tineq}, we can derive upper bounds for the distance variance. 
	%in terms of the variance of the components $X^{(1)},\ldots,X^{(p)}$ and the Gini mean difference of the vector $X$.
	
%	\texcolor{red} {
	\begin{theorem}{\label{th:dimpineq}}
		Let $X \in \R^p$ be a random vector with $\E\|X\|< \infty$ and let $\bX= (X_1,\ldots,X_n)$ denote an i.i.d. sample drawn from $X$. Then
		$$
		\V_n^2(\bX) \leq \frac{1}{n^4} \Big( \sum_{i=1}^n \sum_{j=1}^n \|X_i-X_j\| \Big)^2.
		$$		
Moreover, denoting by $X'$ an independent copy of $X$, we obtain
		$
		\V^2(X) \leq (\E\|X-X'\|)^2.
		$
		
		Further, if $\E \|X\|^2 < \infty$ then 
		$
		\V^2(X) \leq \tr(\Sigma_X),
		$
		where $\Sigma_X$ is the covariance matrix of $X$. 
	\end{theorem}
%}
	
	In the one-dimensional case, Theorem \ref{th:dimpineq} implies that the distance variance is bounded above by the variance and the squared Gini mean difference.  
	
	\begin{corollary}\label{cor:vargmdineq}
		Let $X$ be a scalar random variable with $\E(|X|) < \infty$.  Then, $\V^2(X) \leq \Delta^2(X)$. Moreover, if $\E(|X|^2) < \infty$ then $\V^2(X) \leq \sigma^2(X)$.
	\end{corollary}

	We also note that for continuous variables $X \in \R$, the inequality $T_2(X) \leq T_1(X)$ can be sharpened.
	
	\begin{proposition}\label{prop:T2T1ineq}
		Let $X$ be a real-valued continuous random variable with $\E(|X|^2) < \infty$.  Then, 
		$T_2(X) \leq \tfrac{2}{3} \, T_1(X)$.
	\end{proposition}	
	
	%\medskip
	
	Interestingly, Gini's mean difference and the distance standard deviation coincide for distributions whose mass is concentrated on two points.
	
	%\medskip	
	
	\begin{theorem}\label{th:bernoulli}
		Let $X$ be Bernoulli distributed with parameter $p$. Then
		$$
		\V^2(X) = \Delta^2(X) = 4 \, p^2 \, (1-p)^2.
		$$
		Conversely, if $X$ is a non-trivial random variable for which $\V^2(X) = \Delta^2(X)$ then the distribution of $X$ is concentrated on two points.
	\end{theorem}
	
	%\medskip
	
	For the Bernoulli distribution with $p=\tfrac{1}{2}$, Theorem \ref{th:bernoulli} implies immediately that $\V^2(X)$, $\sigma^2(X)$, and $\Delta^2(X)$ attain the same value, namely, $\tfrac{1}{4}$.  Hence, applying Corollary \ref{cor:vargmdineq} and the dilation property $\V(b X) = |b| \V(X)$ in (C2), we obtain
	
	\begin{corollary}{\label{cor:dvarmax}}
		Let $\mathcal{X}$ denote the set of all real-valued random variables and let $c>0$. Then
		$$
		\max_{X \in \mathcal{X}} \{\V^2(X): \sigma^2(X)=c\} = \max_{X \in \mathcal{X}}\{\V^2(X): \Delta^2(X) = c\} = c,
		$$
		and both maxima are attained by $Z = 2 \, c^{1/2} \, Y$, where $Y$ is Bernoulli distributed with parameter $p = \tfrac12$.
	\end{corollary}	

\noindent	
	This result answers a question raised by G\'{a}bor Sz\'ekely (private communication, November 23, 2015).
	
	We remark that the second implication of Theorem \ref{th:dimpineq} and Theorem \ref{th:bernoulli} also follow from a result for the generalized distance variance in \citep[Proposition 2.3]{lyons2013}. However, our presentation provides a more direct approach to these findings.

	Since distance standard deviation terms appear in the denominator of the distance correlation coefficient, the inequalities derived in this section lead to new properties for the distance correlation. As an example we now state a result, on the behavior of the empirical distance correlation in high dimensions, that can be derived using Theorem \ref{th:dimpineq} (see Appendix \ref{appsec:proofs} for full details). In \citep[Appendix A.1.]{szekely2013}, it is shown under certain assumptions that $\V_n^2(\bX,\bY)$ converges to $1$ almost surely when the dimensions of $\bX$ and $\bY$ tend to infinity, while the sample size $n$ is fixed. We now show that a similar property can be derived when only the dimension of $\bX$ tends to infinity.

	 Before stating the result, we note (see \cite[p. 2776, Eq. (2.18)]{szekely2007}) that the squared standard empirical distance covariance $\V_n^2(\bX,\bY)$ is expressible as
		\begin{equation}
		\begin{aligned}
		\V_n^2(\bX,\bY) &=  \frac{1}{n^2} \sum_{i,j=1}^n \|X_i-X_j\| \|Y_i-Y_j\| \\
		& \qquad + \frac{1}{n^4} \sum_{i,j=1}^n \|X_i-X_j\| \sum_{i,j=1}^n \|Y_i-Y_j\| \\ 
		& \qquad - \frac{2}{n^3} \sum_{i,j,k=1}^n \|X_i-X_j\| \|Y_i-Y_k\|. \label{eq:sampledcov}
		\end{aligned}
		\end{equation}
	The squared standard empirical distance correlation (see \cite[p. 2774, Definition 5]{szekely2007}) is defined as
		\begin{equation}
	\mathcal{R}_n^2(\bX,\bY) =  \frac{\V_n^2(\bX,\bY)}{\V_n(\bX) \V_n(\bY) },\label{eq:sampledcor}
		\end{equation}
		if both $\V_n(\bX)$ and $\V_n(\bY)$ are different from $0$, and $\mathcal{R}_n^2(\bX,\bY)=0$ otherwise. 
		
	\begin{proposition}{\label{prop:samplelim}}
			For fixed $q$, let $Y$ be a $q$-dimensional random vector.  For each $p \in \N$, let $X  = (X^{(1)},\ldots,X^{(p)})^t$ be a $p$-dimensional random vector with $\E \|X\|^2 < \infty$ and i.i.d. coordinates $X^{(1)},\ldots,X^{(p)}$. For fixed $n \in \N$, let $(\bX,\bY) =((X_1,Y_1),\ldots,(X_n,Y_n))$  denote a sample of size $n$ drawn from $(X,Y)$.  Then, almost surely, 
			\begin{equation} \label{eq:dcovhighdim1}
			\lim_{p \to \infty} \frac{{\V}_n^2(\bX,\bY)}{\E\|X-X'\|}  = n^{-3} \,  \sum_{i,j=1}^{n} \|Y_i-Y_j\|
			\end{equation}
			and		
			\begin{equation} \label{eq:dcorhighdim1}
			\lim_{p \to \infty} {{\mathcal{R}}}_n^2(\bX,\bY) = (n-1)^{-1/2} \, \frac{n^{-2}\, \sum_{i,j=1}^{n} \|Y_i-Y_j\|}{{\V}_n(\bY,\bY)}\geq (n-1)^{-1/2}.	
			\end{equation}
	\end{proposition}

To demonstrate the relevance of Proposition \ref{prop:samplelim}, we generate i.i.d. samples $(\bX^{(k)}, \bY^{(k)}) = ((X^{(k)}_1,Y^{(k)}_1),\ldots,(X^{(k)}_n,Y^{(k)}_n))$ of size $n=50$ drawn from $(X,Y)$, where $X \in \R^p$ with $p=100$, $Y \in \R$ and $(X,Y)$ follows a $p+1$-dimensional standard normal distribution with identity covariance matrix. For the average standard distance correlation over $K=10,000$ simulation runs we then obtain $K^{-1} \sum_{k=1}^K \mathcal{R}_n(\bX^{(k)}, \bY^{(k)})= 0.4823$. Considering that we simulated $\bX$ and $\bY$ to be independent, this reveals a heavy bias of the standard distance correlation in this setting, showing that this coefficient is hard to interpret when $p$ is high. Even more, we note that the limiting value of ${{\mathcal{R}}}_n^2(\bX,\bY)$ 
% $(n-1)^{-1/2} \, \frac{n^{-2}\, \sum_{i,j=1}^{n} \|Y_i^{(k)}-Y^{(k)}_j\|}{{\V}_n(\bY^{(k)},\bY^{(k)})}$ 
depends only on the distribution of $Y$ and not on the dependence between $X$ and $Y$. Hence, we can expect similar results for random variables following the same distribution as $Y$ even when they are strongly associated with $X$. Indeed, let ${\bf Z}^{(k)} = (Z^{(k)}_1,\ldots,Z^{(k)}_n)$, where $Z^{(k)}_i = p^{-1/2} \boldsymbol{1}_p' X^{(k)}_i$ with $\boldsymbol{1}_p = (1,\ldots,1)' \in \R^p$. Obviously $Z^{(k)}_i$ shows the same variance as $Y^{(k)}_i$, but now $Z^{(k)}_i$ and $X^{(k)}_i$ are collinear. Yet, we obtain $K^{-1} \sum_{k=1}^K \mathcal{R}_n(\bX^{(k)}, {\bf Z}^{(k)})= 0.5112$, showing only a slight difference to the result in the independent case. For an interpretable version of distance correlation when $p$ is high, we propose to use
	\begin{equation}\label{eq:bcdcor}
	\widehat{\mathcal{R}}_n(\bX,\bY) = \mbox{sign}\left(  \frac{\widehat{\Omega}_n(\bX,\bY)}{\widehat{\V}_n(\bX) \widehat{\V}_n(\bY)} \right)  \, \sqrt{\left| \frac{\widehat{\Omega}_n(\bX,\bY)}{\widehat{\V}_n(\bX) \widehat{\V}_n(\bY)} \right|},
	\end{equation}	
where $\mbox{sign}(t)$ denotes the sign of $t \in \R$ and $\widehat{\Omega}_n(\bX,\bY)$ is defined in \eqref{eq:ugendcov}. 

Since this version is based on the U-statistic estimates of the squared distance covariance and the distance variance, it may be conjectured that it will generally not show a strong bias. Notably, $K^{-1} \sum_{k=1}^K \widehat{\mathcal{R}}_n(\bX^{(k)}, \bY^{(k)})= -0.00511$ and $K^{-1} \sum_{k=1}^K \widehat{\mathcal{R}}_n(\bX^{(k)}, {\bf Z}^{(k)})= 0.2897$; the population versions can be explicitly calculated using Corollary 3.2 and Corollary 3.3 in \cite{Dueck2014} and are given by $\mathcal{R}(X,Y) = 0$ and $\mathcal{R}(X,Z) \approx 0.2987$.

Examples in which the dimension of $X$ is $100$ or larger and $Y$ is univariate occur in the analysis of genetic data, where it is often the goal to assess the association of a large number of molecular markers with some univariate clinical response, such as the development of a certain disease or response to treatment. One common approach \cite{goeman2007} for this kind of data is to test for the association of the response with interesting {\it sets of markers} which may for example be defined  {\it via} gene pathways or gene ontology (GO) \cite{GO} terms. While hypothesis testing itself gives little information about the effect size, distance correlation offers a way to quantify the strength of association between sets and univariate responses. Proposition \ref{prop:samplelim} and the above considerations yield that the bias-corrected estimate (\ref{eq:bcdcor}) is to be preferred over the standard estimator (\ref{eq:sampledcor}) in these situations.

\section{Properties of the distance standard deviation in one dimension}
\label{sec:properties}		
	
	%	\section{New representations for the distance variance}{\label{sec:altrep}}
	
	The representation of $\V^2$ given in (\ref{eq:dvartwo}), although more applicable than the expression given in equation (\ref{eq:dvar}), is undefined for random vectors with infinite second moments.  This problem can be circumvented by considering the representation 
	\begin{equation} \label{eq:dvarthree}
	\V^2(X) = \Delta^2(X) + W(X),
	\end{equation}
	where 
	$$
	W(X) =  \E \Big[\|X-X'\| \cdot \big(\|X-X'\| - 2 \, \|X-X''\| \big) \Big].
	$$
	Note that since $0 \le \V^2(X) \le \Delta^2(X)$ then $W(X) \le 0$ and $|W(X)| \le \Delta^2(X)$; since $\Delta(X)$ exists under the assumption of finite first-order moments, then so does $W(X)$.

	In the one-dimensional case, (\ref{eq:dvarthree}) gives rise to other representations that lead to crucial results about the distance standard deviation.
	
	\begin{theorem}\label{th:repdiff}
		Let $X$ be a real-valued random variable with $\E|X| < \infty$, and let $X$, $X'$, $X''$, and $X'''$ be i.i.d. %Then $\V^2(X)$, the squared distance standard deviation, can be expressed as follows:
		\begin{itemize}
			\item[(i)]  Let $X_{1:4} \leq X_{2:4} \leq X_{3:4} \leq X_{4:4} $ be the order statistics of the quadruple $(X,X',X'',X''')$.  Then, 
		\begin{equation}
		\V^2(X) = \frac{2}{3} \E [(X_{3:4}-X_{2:4})^2]. \label{eq:dvarrepcum2}
		\end{equation}	
			\item[(ii)]  Let $F$ be the cumulative distribution function of $X$.  Then,
		\begin{equation}
		\V^2(X)= 8 \operatornamewithlimits\iint\limits_{-\infty<x<y<\infty}  F^2(x) (1-F(y))^2 \dd x \, \dd y. \label{eq:dvarrepcum}
		\end{equation}
			\item[(iii)]  Let $t_+ = \max(t,0)$, $t \in \R$.  Then, 
			\begin{equation}
			\V^2(X) = \Delta^2(X)- 8 \, \E[(X-X')_{+}\,(X''-X)_+]. \label{eq:dvarfour2} 
			\end{equation}	
			\item[(iv)]  Let $X_{1:3} \leq X_{2:3} \leq X_{3:3}$ be the order statistics of the triple $(X,X',X'')$.  Then,
			\begin{equation}
			\V^2(X) = \Delta^2(X)- \tfrac43 \, \E[(X_{2:3}-X_{1:3})\,(X_{3:3}-X_{2:3})]. \label{eq:dvarfour}
			\end{equation}
		\end{itemize}
	\end{theorem}
	
Important properties of $\V$ following from equation \eqref{eq:dvarrepcum2} are discussed in Theorem \ref{th:disp}, and motivation for the representations provided in \eqref{eq:dvarrepcum}-\eqref{eq:dvarfour} are given in the supplementary material.

\begin{theorem}\label{th:disp}	
	The functional $\V$ is an axiomatic measure of spread, i.e.,
		\begin{itemize}
		\item[(C1)] $\V(X) \geq 0$,
		\item[(C2)] $\V(a+bX) = |b| \, \V (X)$ for all $a,b \in \R$, and 
		\item[(C3)] $\V(X) \leq \V(Y)$ if for all $0 < \alpha \leq \beta < 1$,
		\begin{equation*}
		%\label{eq:dispersive}
		F^{-1}(\beta) - F^{-1}(\alpha) \leq  G^{-1}(\beta) - G^{-1}(\alpha),
		\end{equation*}
		where $F$ and $G$ are the cumulative distribution functions of $X$ and $Y$, respectively, and $F^{-1}$ and $G^{-1}$ are the corresponding right-continuous inverses.
	\end{itemize}
\end{theorem}	

	Applying \citep[Theorem 3.B.7]{Shaked1994}, we obtain the following corollary of Theorem \ref{th:disp}.

	\begin{corollary}{\label{cor:logconcave}}
		Let $X$ be a random variable with a log-concave density. Then
		$
		\V(X+Y) \geq \V(X)
		$
		for any random variable $Y$ independent of $X$.
	\end{corollary}
	In particular, if $X$ and $Y$ are independently distributed, continuous, random variables with log-concave densities, then
	\begin{equation} \label{eq:propmax}
	\V(X+Y) \geq \max(\V(X),\V(Y)).
	\end{equation}
	It is well known, both for the standard deviation and for Gini's mean difference, that assertions analogous to (\ref{eq:propmax}) hold without restrictions on the distributions of $X$ and $Y$.

	We now show, however, that this property does not hold generally for the distance standard deviation, $\V$, thereby answering a second question raised by G\'abor Sz\'ekely (private communication, November 23, 2015).  
	
	\begin{example}{\label{ex:counterex1}}
		Let $X$ be Bernoulli distributed with parameter $p=\tfrac12$ and let $Y$ be uniformly distributed on the interval $[0,1]$ and independent of $X$. Then $\V(X) > \V(X+Y)$.  
	\end{example}
	
%	\smallskip
	
	Other common properties of the classical standard deviation and Gini's mean difference concern differences and sums of independent random variables. Notably, it is well-known that,
	$\Delta(X+Y) = \Delta(X-Y)$ and $\sigma(X+Y) = \sigma(X-Y)$
	for any independent random variables $X$ and $Y$ for which these expressions exist.  
	On the other hand, these properties do not hold in general for the distance standard deviation.
	
	\begin{example}{\label{ex:nonsym}}
		Let $X$ and $Y$ be independently Bernoulli distributed with parameter $p \neq \tfrac12$. Then $\V(X+Y) > \V(X-Y)$.
	\end{example}

	While $\Delta(X)$, $\sigma(X)$ and $\V(X)$ are all measures of spread in the sense of \cite{Bickel2012}, Examples \ref{ex:counterex1} and \ref{ex:nonsym} and the comparison of the asymptotic relative efficiencies in Section \ref{sec:estimator} suggest that there are substantial differences between these coefficients as measures of spread. To provide further understanding of these differences, we now derive representations that enable graphical comparisons of these three measures.
	
	For this purpose, we apply equation (\ref{eq:dvarrepcum2}) to derive a new empirical version for distance variance which is distinct from $\V_n^2(\bX)$ and $\widehat{\V}_n^2(\bX)$, as follows. For an i.i.d. sample $X_1,\ldots,X_n$ of real-valued random variables, denote by $D_{i:n} = X_{(i+1):n} - X_{i:n}$, $i=1,\ldots,n-1$ the $i$th {\it spacing of} $\bX = (X_1,\ldots,X_n)$.

	\begin{proposition}\label{cor:newsampledistvar}
	Let $X$ be a real-valued random variable with $\E(|X|) < \infty$ and let $\bX=(X_1,\ldots,X_n)$ be an i.i.d. sample from $X$.  Then, a strongly consistent empirical version for $\V^2(X)$ is 
	\begin{equation}
	\label{eq:newsampledvar}
	\U_n^2(\bX) = \binom{n}{2}^{-2} \sum_{i,j=1}^{n-1} \big(\min(i,j)\big)^2  \big(n-\max(i,j)\big)^2 D_{i:n} D_{j:n}.
	\end{equation}
	%		where $D_{k:n}= X_{{k+1}:n} - X_{k:n}$ denotes the $k$th sample spacing of $X$, $1 \le k \le n-1$.
\end{proposition}

 Let $D = (D_{1:n},\ldots,D_{(n-1):n})$ denote the vector of spacings, then we can write the quadratic form in (\ref{eq:newsampledvar}) as
	$
	\U_n^2(\bX) = D^t \, V \, D,
	$
	where the $(i,j)$th element of the matrix $V$ is 
	$$
	V_{i,j} = \binom{n}{2}^{-2} \, \big(\min(i,j)\big)^2 \,  \big(n-\max(i,j)\big)^2.
	$$	
 Both the squared empirical Gini mean difference and the empirical variance (see equations (\ref{eq:empsd}) and (\ref{eq:empgini}))
	can also be expressed as quadratic forms in the spacings vector $D$; specifically, 
		$
	\widehat{\Delta}_n^2(\bX) =  D^t \, G \, D$ and  $\widehat{\sigma}_n^2(\bX) =  D^t \, S \, D$,
		where the elements of $G$ and $S$ are given by
	$$
	G_{i,j} = \binom{n}{2}^{-2} \, i\, j  \, (n-i) \  (n-j)
	$$
	and
	$$
	S_{i,j} = \frac{1}{2} \binom{n}{2}^{-1} \, \min(i,j) \, \big(n-\max(i,j)\big).
	$$

	\begin{figure}[t!]
		\centering
		\captionsetup{width=0.8\textwidth}
		\includegraphics[width=1.0\textwidth]{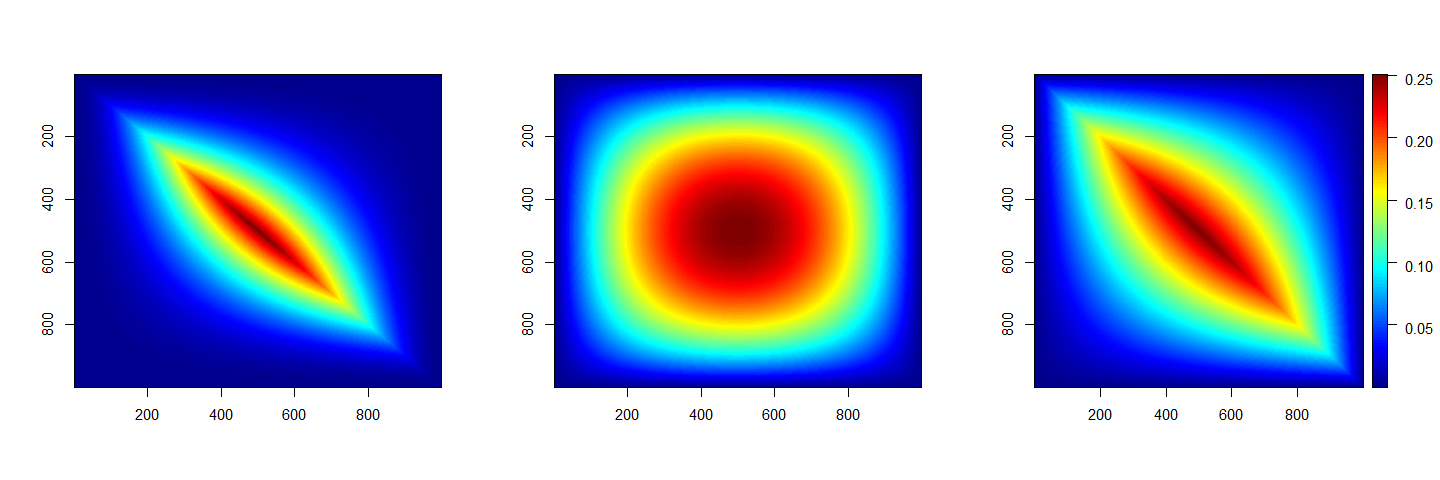}
		\caption{\small Illustration of, from left to right, the empirical distance variance $\U^2_n$, the squared empirical Gini mean difference $\widehat{\Delta}_n^2$, and the empirical variance $\widehat{\sigma}_n^2$ {\it via} their respective quadratic form matrices $V$, $G$, and $S$ for sample size $n=1,000$. The coordinate $(i,j)$ corresponds to the $(i,j)$th entry of the corresponding matrix, and the size of the corresponding matrix element is specified {\it via} color code.}
		\label{figure1}
	\end{figure}

Comparing $\U^2_n$, $\widehat{\Delta}_n^2$, and $\widehat{\sigma}_n^2$ clearly is equivalent to comparing the matrices $V$, $G$ and $S$. We use this fact to graphically illustrate differing features of $\V$, $\Delta$, and $\sigma$ by plotting the values of the underlying matrices; see Figure \ref{figure1}. These plots provide a descriptive explanation as to why $\V$ and $\Delta$ are more suitable for heavy-tailed distribution since they place smaller weight than $\sigma$ on extreme spacings; in particular, $V_{n-1,n-1} = G_{n-1,n-1} = 4n^{-2}$ while $S_{n-1,n-1} = n^{-1}$.  
	
	On the other hand, the {\it shape} of the plot for $\V$ resembles the plot for $\sigma$ more than the plot for $\Delta$. Both $\V$ and $\sigma$ place comparably high weights on the diagonals; specifically, the highest entries in row $i$ are the diagonal entries $V_{i,i}$ and $S_{i,i}$, respectively, and the highest entry in row $i$ for Gini's mean difference is $G_{i,\lfloor n/2\rfloor}$. The intrinsic reason for this finding is that $\V^2$ and $\sigma^2$ are sums of {\it quadratic} differences of the observations, while $\Delta$ is a sum of the untransformed differences.  
	
	Since $\V$ and $\sigma$ place comparably high weights on the diagonals, one may conjecture that these measures are sensitive to variability in the spacings.  Indeed, Yitzthaki \citep[p. 291]{yitzhaki2003} points out that ``the more equal are the distances between adjacent observations, the lower the variance''.  This statement also holds for the distance standard deviations of light-tailed distributions. 
	
	For example, consider a distribution which is concentrated on three points, $\{0,\delta,1\}$ and attains each value with probability $\tfrac13$, similar to \citep[p. 291]{yitzhaki2003}. For this distribution, $\Delta = \tfrac{4}{9}$, irrespective of the value of $\delta$. On the other hand, \:$\V^2 = 4 \,(\tfrac13 + \delta^2+(1-\delta)^2)/27$ \, and \: $\sigma^2 = 2\,(1 + \delta^2+(1-\delta)^2)/9$. Hence, both measures are maximal for $\delta \in \{0,1\}$ and minimal for $\delta = \tfrac12$, i.e. when the difference between the spacings is maximal and minimal, respectively.
		
	The sensitivity of $\V$ to differences between the spacings extends to other light-tailed distributions. For example, the most broadly spread distribution on $[0,1]$ arguably is $B(1,\tfrac12)$, the Bernoulli distribution with $p=\tfrac12$. The ratio of the population values between $B(1,\tfrac12)$ and the uniform distribution on $[0,1]$ are $3/2 = 1.5$ ($\Delta$), $\sqrt{3} \approx 1.73$ ($\sigma$) and $\sqrt{45/8} \approx 2.37$ ($\V$); see Appendix \ref{sec:distr} for the value of $\V$ for the uniform distribution. 
	
	In the case of $\sigma$, Yitzthaki \citep[p. 291]{yitzhaki2003} states that the sensitivity to differences between the spacings ``is translated to sensitivity to extreme observations.'' However, this is not the case for $\V$ as it places lower weights than $\Delta$ on the outer spacings, resulting in less sensitive behavior to extreme observations. For a demonstrative example, consider the $t_\nu$-distributions. Specifically, the ratio of the respective population values between the $t_3$- and $t_5$-distribution is \, $(54 \sqrt{3})/(35 \sqrt{5}) \approx 1.20$ \, for $\Delta$ and $\sqrt{9/5} \approx 1.34$ for $\sigma$; for $\V$, numerical evaluation yields that the ratio is approximately $1.12$.

To summarize, the distance standard deviation is very sensitive to variability in the central spacings of an observation but relatively insensitive to changes in the extreme spacings.  Consequently, the behavior of $\V$ resembles the behavior of $\sigma$ for light-tailed distributions where the central spacings are relatively large compared to the extreme spacings. On the other hand, the low sensitivity of $\V$ to extreme observations makes it a very good measure of spread for heavy-tailed distributions.

\section{Discussion} \label{sec:discussion}

In this work, we have studied the statistical properties of the distance standard deviation, which arises as a special case of the distance covariance introduced by Sz\'{e}kely, et al. \cite{szekely2007}. Notably, we have demonstrated that the empirical distance standard deviation has appealing statistical properties: it is less vulnerable to outliers and generally more appropriate for heavy-tailed distributions, more so than classical alternatives such as the mean deviation and Gini's mean difference.

The distance standard deviation, and the distance correlation coefficient, may be seen as natural statistics for analyzing multivariate distributions, providing an alternative to the classical second-moment statistics and also being potentially more appropriate in light of their statistical properties.

For multivariate random variables, the distance standard deviation summarizes the spread as a single value. While this can be useful for many applications (see, e.g., the example on multivariate statistical quality control in Section \ref{sec:app}), a referee has noted that the covariance matrix provides richer information, such as the spread of single components and the association between different components.  To obtain an analogue of the covariance matrix based on the concept of distance covariance, one can define the {\it distance covariance matrix}, 
	$$
		\Gamma_X =	\left(\V^2(X_i,X_j)\right)_{i,j=1,\ldots,p} \, ,
	$$
where $X_1,\ldots,X_p$ are the components of $X$; a related concept is that of distance multivariance \cite{Boettcher2019}, which allows for testing the mutual independence of more than two sets of random vectors. It can be shown that $\Gamma_X$ is positive semidefinite. Investigating properties of the distance covariance matrix is a promising direction for further research.

\section*{Acknowledgments}
We are grateful to the reviewers, the associate editor and the editor for constructive and insightful comments that led to numerous revisions in the manuscript.  Dominic Edelmann is funded by the German Research Foundation (Project No. 417754611).

	\newpage
	
\begin{frontmatter}
	
	\title{Supplementary Material to \\ ``The Distance Standard Deviation''}
	\runtitle{Supplementary Material}
	
	\begin{aug}
		\author{\fnms{Dominic} \snm{Edelmann}\thanksref{t1},\ead[label=e1]{dominic.edelmann@dkfz-heidelberg.de}}
		\author{\fnms{Donald} \snm{Richards},\ead[label=e2]{richards@stat.psu.edu}}
		\and
		\author{\fnms{Daniel} \snm{Vogel}\ead[label=e3]{daniel.vogel@abdn.ac.uk}}
		
		\thankstext{t1}{Corresponding author}
		\runauthor{D. Edelmann, D. Richards and D. Vogel}
		
		\affiliation{German Cancer Research Center, Pennsylvania State University, and University of Aberdeen}
		
		\address{German Cancer Research Center\\ %(Deutsches Krebsforschungszentrum), 
			Im Neuenheimer Feld 280\\
			69120 Heidelberg \\
			Germany\\
			\printead{e1}\\}
		
		\address{Department of Statistics\\ 
			Pennsylvania State University\\ 
			University Park, PA 16802\\ 
			U.S.A.\\
			\printead{e2}\\}
		
		\address{
			Institute for Complex Systems and Mathematical Biology\\ University of Aberdeen\\ Aberdeen AB24 3UE\\ U.K. \\
			\printead{e3}\\
		}
		
	\end{aug}

	\begin{abstract}
		This supplementary material consists of six appendices. Appendix \ref{appsec:proofs} contains various proofs to results in the main paper. Appendix \ref{appsec:table1} contains details on the derivation of Table \ref{tab:estimators.asymptotic} in the main document. In Appendix \ref{app:clt}, a limit theorem for the squared distance covariance is stated, under weaker assumptions than known previously. Simulation results for permutation-based two-sample scale tests are provided in Appendix \ref{app:permtest}. Appendix \ref{app:furtherprop} gives additional theoretical results for the distance standard deviation in one dimension. Appendix \ref{sec:distr} tabulates the distance variances for a collection of distributions. Equations, theorems, tables, etc., are referred to as in the main document, e.g., (\ref{eq:dcov}) refers to Equation (\ref{eq:dcov}) in the main document, while  (\ref{eq:kernelexp}) refers to Equation (\ref{eq:kernelexp}) in Appendix A. 
	\end{abstract}

\end{frontmatter}

%%%%%%%%%%%
%
%
%
%
%
%
%
%
%

\appendix

\section{Proofs}
\label{appsec:proofs}

\begin{proof}[Proof of Lemma \ref{lem:kernel}]
	Expanding the kernel function given in (\ref{eq:kernel}) yields
	\begin{align}
	h(X_1,X_2,X_3,X_4)&=\frac{1}{4} \, \sum_{\substack{1 \leq i,j \leq 4\\ i \neq j}}  \|X_i - X_j \|^2 \nonumber \\ 
	&\quad - \frac{1}{4} \, \sum_{\substack{1 \leq i,j,k \leq 4 \\ i \neq j, i \neq k}}  \|X_i - X_j \| \, \|X_i - X_k \| \label{eq:kernelexp}\\
	&\quad + \frac{1}{24}  \sum_{\substack{1 \leq i,j,k,l \leq 4\\ i \neq j, k \neq l}} \|X_i - X_j \| \,  \|X_k - X_l \|. \nonumber
	\end{align}
	Next, we will expand each term in the representation (\ref{eq:kernelexp}).  Denote by $A \dot\cup B$ the disjoint union of two sets $A$ and $B$, and also denote by $\#(A)$ the cardinality of $A$.  By partitioning the set $\big\{1 \leq i,j,k,l \leq 4, i \neq j, l \neq k \big\}$ into the subsets in which either zero, one, or two of the elements $(i,j)$ and $(k,l)$ coincide, we obtain 
	\begin{align*}
	\big\{1 \leq i,j,k,l & \leq 4, i \neq j, l \neq k \big\} \\
	&= \big\{1 \leq i,j,k,l \leq 4\, \big| \,  i,j,k,l \text{ all different}\big\} \\  
	&\qquad  \dot \cup \big\{1 \leq i,j,k,l \leq 4\, \big| \,  i \neq j, k \neq l, \, \#\big(\{i,j\} \cap \{k,l\}\big) = 1 \big\} \\ \
	&\qquad \dot \cup \big\{1 \leq i,j,k,l \leq 4 \, \big| \, i \neq j, k \neq l, \, \#\big(\{i,j\} \cap \{k,l\}\big) = 2 \big\}.
	\end{align*}
	By using the notation $a_{ij} := \|X_i - X_j \|$, $i,j \in \{1,\ldots,n\}$ and applying a symmetry argument, the third summand in the representation (\ref{eq:kernelexp}) can be written in the form, 
	\iffalse
	\begin{align*}
	\sum_{\substack{1 \leq i,j,k,l \leq 4\\ i \neq j, k \neq l}} a_{ij} a_{kl} &= \sum_{\substack{1 \leq i,j,k,l \leq 4\\ i,j,k,l \text{ all different}}} a_{ij} a_{kl} \\
	&\qquad + \sum_{\substack{1 \leq i,j,k \leq 4\\ i,j,k \text{ all different}}} (a_{ij} a_{ik} + a_{ij} a_{ik} + a_{ij} a_{jk} + a_{ij} a_{jk}) \\
	&\qquad + \sum_{\substack{1 \leq i,j \leq 4\\ i \neq j}} (a_{ij} a_{ij} + a_{ij} a_{ij}),
	\end{align*}
	and by applying a symmetry argument, the latter equation reduces to 
	\fi
	\begin{multline}
	\sum_{\substack{1 \leq i,j,k,l \leq 4\\ i \neq j, k \neq l}} a_{ij} a_{kl} \\ 
	= \sum_{\substack{1 \leq i,j,k,l \leq 4\\ i,j,k,l \text{ all different}}} a_{ij} a_{kl} 
	+  4 \, \sum_{\substack{1 \leq i,j,k \leq 4\\ i,j,k \text{ all different}}} a_{ij} a_{ik}  
	+ 2 \sum_{\substack{1 \leq i,j \leq 4\\ i \neq j}} a_{ij}^2. \label{eq:kernthirdsummand}
	\end{multline}	
	Next, the set decomposition,  
	\begin{align*}
	\big\{1 \leq i,j,k\leq 4, i \neq j, i \neq k\big\}  &= \big\{1 \leq i,j,k\leq 4 \, \big|  \, i,j,k \text{ all different}\big\} \\ 
	& \qquad \dot \cup\big\{1 \leq i,j,k\leq 4 \, \big| \, i \neq j, i \neq k, j= k\big\},
	\end{align*}
	yields for the second summand in (\ref{eq:kernelexp}) the expression, 
	\begin{align}
	&\sum_{\substack{1 \leq i,j,k \leq 4\\ i \neq j, i \neq k}} a_{ij} a_{ik}= \sum_{\substack{1 \leq i,j,k \leq 4\\i,j,k \text{ all different}}} a_{ij} a_{ik} + \sum_{\substack{1 \leq i,j \leq 4\\ i \neq j}} a_{ij}^2. \label{eq:kernsecsummand}
	\end{align}
	Inserting  (\ref{eq:kernthirdsummand}) and (\ref{eq:kernsecsummand}) into (\ref{eq:kernelexp}) completes the proof.
\end{proof}

\begin{proof}[Proof of Theorem \ref{th:asymptotics}]
	By \citep{Huo2015}, it follows that 
	\begin{equation} \label{eq:kernpos}
	h(X_1,X_2,X_3,X_4) =	\frac{1}{4} \sum_{i,j=1}^4  \widetilde{A}_{ij}^2,
	\end{equation}
	where 
	\begin{equation} \label{eq:quadratic}
	\widetilde{A}_{ij} = \begin{cases} 
	a_{ij}-\frac{1}{2} \sum_{i=1}^n a_{ij} - \frac{1}{2} \sum_{j=1}^n a_{ij} + \frac{1}{6} \sum_{i,j=1}^n a_{ij}, \quad & i \neq j, \\ 
	0, & i = j. 
	\end{cases}
	\end{equation}
	%	where $a_{ij} = \|X_i - X_j\|$.
	By representation (\ref{eq:kernpos}), we see that $h(X_1,X_2,X_3,X_4) \geq 0$. Moreover, an elementary computation yields
	$$
	\sum_{\substack{1 \leq i,j \leq 4\\ i \neq j}}^4  \|X_i - X_j \|^2 - \sum_{\substack{1 \leq i,j,k \leq 4\\i,j,k \text{ all different}}} \|X_i - X_j \| \, \|X_i - X_k \| \leq 0,
	$$
	hence 
	$$
	h^2(X_1,X_2,X_3,X_4) \leq \frac{1}{24^2} \Bigg( \sum_{\substack{1 \leq i,j,k,l \leq 4\\ i,j,k,l \text{ all different}}} \|X_i - X_j \| \,  \|X_k - X_l \| \Bigg)^2.
	$$
	By the Cauchy-Schwarz inequality, $(\sum_{i=1}^n a_i)^2  \leq n \sum_{i=1}^n a_i^2$ for $a_1,\ldots,a_n \in \R$, hence
	$$
	h^2(X_1,X_2,X_3,X_4) \leq \frac{1}{24}  \sum_{\substack{1 \leq i,j,k,l \leq 4\\ i,j,k,l \text{ all different}}} \|X_i - X_j \|^2 \,  \|X_k - X_l \|^2;
	$$	
	therefore, 
	\begin{align}
	\E h^2(X_1,X_2,X_3,X_4) &\leq \E |X_1 - X_2|^2 \, |X_3 - X_4|^2  \\
	&\leq C \E |X_1|^2 \, \E |X_2|^2 \\
	&= C \left(\E |X_1|^2 \right)^2, \label{eq:h2bounded}
	\end{align}
	where $C$ is a constant and the last inequality follows from the fact that $\E |X_1 -X_2|^2 \leq 2 \E |X_1|^2 + 2 \E |X_2|^2$.
	
	Hence $\E [h^2(X_1,X_2,X_3,X_4)] < \infty$, and we deduce from a classical result of Hoeffding \citep[Theorem 7.1]{Hoeffding1948} that 
	$$
	\sqrt{n}\big(\widehat{\V}_n^2(\bX) - \V^2(X)\big) \xrightarrow{d} N\big(0, 16 \, \E [h_1^2(X)] \big).
	$$
	The asymptotic distribution of the biased version $\V_n^2(\bX)$ can now be obtained from the asymptotic distribution of $\widehat{\V}_n^2(\bX)$ since
	$$
	\V_n^2(\bX) =  \widehat{\V}_n^2(\bX) + O_p (n^{-1}).
	$$
	Finally, (\ref{eq:Vnasymptoticgammahat}) is obtained by applying Slutsky's Theorem to (\ref{Vnasymptoticdistn}).
\end{proof}

\begin{proof}[Proof of Theorem \ref{th:testlim}]
	For proving (i), we note that
	\[
	\widehat\xi_{p}(\bX_n, \bY_m) = \frac{n}{n+m} \widehat\xi(\bX_n) + \frac{m}{n+m} \widehat\xi(\bY_m)
	\cip \frac{r}{1+r} \xi_X + \frac{1}{1+r} \xi_Y.
	\]
	Now abbreviating $\widehat\xi_{p}(\bX_n, \bY_m)$  by $\widehat\xi_{p}$ and defining 
	$$
	\xi_{\tiny \mbox{tot}} = \frac{r}{1+r} \xi_X + \frac{1}{1+r} \xi_Y,
	$$
	we obtain 
	\begin{align*}
	&\widehat{T}_\V  = \sqrt{\frac{n \, m}{n + m}} \, \frac{\widehat\V_n(\bX_n) - \widehat\V_m(\bY\!_m)}{ \sqrt{\widehat{\xi}_p}} \\
	&=\sqrt{\frac{n \, m}{n + m}} \, \frac{\widehat\V_n(\bX_n) - \V(X)}{ \sqrt{\widehat{\xi}_p}} - \sqrt{\frac{n \, m}{n + m}} \, \frac{\widehat\V_n(\bY_m) - \V(X)}{ \sqrt{\widehat{\xi}_p}} \\
	&= \Bigg(\sqrt{\frac{ m}{n + m}} \, \frac{\sqrt{n} \, (\widehat\V_n(\bX_n) - \V(X))}{\sqrt{\xi_{\tiny \mbox{tot}}}} - \sqrt{\frac{n}{n + m}} \, \frac{\sqrt{m} \, (\widehat\V_n(\bY_m) - \V(X))}{ \sqrt{\xi_{\tiny \mbox{tot}}}} \Bigg) \sqrt{\frac{\xi_{\tiny \mbox{tot}}}{\widehat{\xi}_p}}.
	\end{align*}
	
	The first term in brackets can be written as
	\begin{equation} \label{eq:termone}
	\sqrt{\frac{ m}{n + m}} \, \frac{\sqrt{n} \, (\widehat\V_n(\bX_n) - \V(X))}{\sqrt{\xi_X}} \sqrt{\frac{\xi_X}{\xi_{\tiny \mbox{tot}}}}.
	\end{equation}
	Since $\sqrt{m/(n + m)} \to 1/(r+1)$, it follows by Corollary \ref{cor:asymptotics} that (\ref{eq:termone}) converges in distribution to $N(0,\xi_X/(r \xi_X + \xi_Y))$. Using similar arguments, we can show that the second term in brackets converges in distribution to $N(0,r \xi_Y/(r \xi_X + \xi_Y))$. Using the fact that $\bX_n$ and $\bY_m$ are mutually independent and applying Slutsky's theorem concludes the proof.
	
	For proving (ii), we note that if $\V(X) < \V(Y)$, by Corollary \ref{cor:asymptotics}, 
	$$
	\frac{\widehat\V_n(\bX_n) - \widehat\V_m(\bY\!_m)}{ \sqrt{\widehat{\xi}_p}} \cip c < 0.
	$$
	Moreover $m n/(n + m) \to \infty$ for $n/m \to r >0$.  Finally, (iii) follows analogously and this concludes the proof.
\end{proof}

\begin{proof}[Proof of Theorem \ref{th:ARE}]
	We show that  
	
	\smallskip
	\noindent
	(a) $\widehat{T}_1(\bX_n,\bY_m^{(n,m)}) \cid N(-\Lambda s_1(X)/\sqrt{\xi_1}, 1)$ and 
	
	\smallskip
	\noindent
	(b) $\widehat{T}_2(\bX_{[\rho n]}, \bY_{[\rho m]}^{(n,m)}) \cid N(-\Lambda s_1(X)/\sqrt{\xi_1}, 1)$. 
	
	In proving (a), we write $\widehat\xi_{p,1}$ to denote $\widehat\xi_{p,1}(\bX_n, \bY_m^{(n,m)})$ and $s_1$ to denote $s_1(X)$. Observe that
	\[
	\widehat\xi_{p,1} = \frac{n}{n+m} \widehat\xi_1(\bX_n) + \frac{m}{n+m} \lambda_{n,m}^2 \widehat\xi_1(\bZ_m)
	\cip \frac{r}{1+r} \xi_1 + \frac{1}{1+r} \xi_1 = \xi_1
	\]
	under the conditions of the theorem. Also, a simple algebraic calculation verifies that 
	\begin{align*}
	\widehat{T}_1(\bX_n, \bY_m^{(n,m)}) &= \sqrt{\frac{n \,m}{n+m}} \frac{{s_n^{(1)}}(\bX_n) - s_m^{(1)}(\bY_m^{(n,m)})}{\sqrt{\widehat{\xi}_{p,1}}} \\
	%&= \ \sqrt{\frac{n \,m}{n+m}} \frac{s_n^{(1)}(\bX_n) - s_1}{\sqrt{\widehat{\xi}_{p,1}}} \ - \ \sqrt{\frac{n \,m}{n+m}} \frac{s_n^{(1)}(\bY_m^{(n,m)}) - \lambda_{n,m}s_1}{\sqrt{\widehat{\xi}_{p,1}}} \]\[ \qquad \qquad - \ \sqrt{\frac{n \,m}{n+m}} \frac{\lambda_{n,m}s_1 - s_1}{\sqrt{\widehat{\xi}_{p,1}}} \\
	&= \Bigg(\sqrt{\frac{m}{n+m}} \frac{\sqrt{n}(s_n^{(1)}(\bX_n) - s_1)}{\sqrt{\xi_1}} \\
	& \qquad - \ \sqrt{\frac{n}{n+m}} \lambda_{n,m} \frac{\sqrt{m}(s_m^{(1)}(\bZ_m) - s_1)}{\sqrt{\xi_1}} \\
	& \qquad\quad - \ \sqrt{\frac{n \,m}{n+m}} (\lambda_{n,m}-1) \frac{s_1}{\sqrt{\xi_1}} \Bigg) \sqrt{\frac{\xi_1}{\widehat{\xi}_{p,1}}}.
	\end{align*}
	The first term in brackets converges in distribution to $\sqrt{1/(r+1)} N(0,1)$, the second term converges in distribution to $\sqrt{r/(r+1)} N(0,1)$. Since the sequences $(X_i)_{i\in\N}$ and $(Z_i)_{i\in\N}$ are mutually independent, the first two terms converge jointly in distribution to the limit $N(0,1)$. The third term in brackets converges to $-\frac{s_1(X)}{\sqrt{\xi_1}}\Lambda$. Statement (a) now follows by applying Slutsky's theorem.
	
	%
	%
	% part (b)
	To prove assertion (b), we abbreviate $\widehat\xi_{p,2}(\bX_{[\rho n]}, \bY_{[\rho m]}^{(n,m)})$ by $\widehat\xi_{p,2}$ and note that
	\begin{eqnarray*}
		\widehat\xi_{p,2} &=& \frac{[\rho n]}{[\rho n]+[\rho m]} \widehat\xi_{2}(\bX_n) + \frac{[\rho m]}{[\rho n]+[\rho m]} \lambda_{n,m}^2 \widehat\xi_{2}(\bZ_m) \\
		&\cip& \frac{r}{1+r} \xi_2 + \frac{1}{1+r} \xi_2 = \xi_2
	\end{eqnarray*}
	under the conditions of the theorem.  Denoting $s_2(X)$ by $s_2$, and proceeding analogously to (a), we obtain
	\begin{align*}
	\widehat{T}_2(\bX_{[\rho n]},\bY_{[\rho m]}^{(n,m)}) &=
	\sqrt{\frac{[\rho n] \,[\rho m]}{[\rho n]+[\rho m]}} \frac{s_{[\rho n]}^{(2)}(\bX_{[\rho n]}) - s_{[\rho m]}^{(2)}(\bY_{[\rho m]}^{(n,m)})}{\sqrt{\widehat\xi_{p,2}}} \\
	&= \Bigg(\sqrt{\frac{[\rho n] \,[\rho m]}{[\rho n]+[\rho m]}} \frac{s_{[\rho n]}^{(2)}(\bX_{[\rho n]}) - s_2}{\sqrt{\xi_{2}}} \\
	& \qquad - \ \sqrt{\frac{[\rho n] }{[\rho n]+[\rho m]}} \lambda_{n,m} \frac  {\sqrt{[\rho m]}(s_{[\rho m]}^{(2)}(\bZ\!_{[\rho m]}) - s_2)}{\sqrt{\xi_{2}}} \\
	& \qquad \quad - \ \sqrt{\frac{[\rho n] \,[\rho m]}{[\rho n]+[\rho m]}} (\lambda_{n,m}-1) \frac{s_2}{\sqrt{\xi_{2}}}\Bigg) \sqrt{\frac{\xi_2}{\widehat{\xi}_{p,2}}}.
	\end{align*}
	By the same reasoning as in (a), the sum of the first two terms in brackets converges in distribution to $N(0,1)$, and the third term converges in probability to 
	\[
	-\sqrt{\rho}\, \Lambda \, \frac{s_2(X)}{\sqrt{\xi_2}} =  - \frac{\sqrt{\xi_2} s_1(X)}{\sqrt{\xi_1} s_2(X)} \Lambda \frac{s_2(X)}{\sqrt{\xi_2}}
	= - \Lambda \frac{s_1(X)}{\sqrt{\xi_1}}.
	\]
	Applying Slutsky's theorem completes the proof.
\end{proof}

\begin{proof}[Proof of Lemma \ref{lem:Tineq}]
	\Pro 
	First note that
	\begin{align*}
	T_{3,n}(\bX) &= \frac{1}{n^3}  \sum_{i=1}^n \sum_{j=1}^n \sum_{k=1}^n  \|X_i-X_j\| \cdot \|X_i-X_k\| \\  &= \frac{1}{n^3} \sum_{i=1}^n \Big( \sum_{j=1}^n \|X_i-X_j\| \Big) ^2.
	\end{align*}
	%Since $(\sum_{i=1}^n a_i)^2  \leq n \sum_{i=1}^n a_i^2$ for any nonnegative numbers $a_1,\ldots,a_n$, an inequality which can be easily proved by induction, it follows that 
	Applying the Cauchy-Schwarz inequality
	%, $(\sum_{i=1}^n a_i)^2  \leq n \sum_{i=1}^n a_i^2$ for all $a_1,\ldots,a_n \allowbreak \in \R$; applying this inequality 
	to the sums that define $T_{1,n}$, $T_{2,n}$, and $T_{3,n}$, we obtain 
	\begin{align*}
	T_{2,n}(\bX)&=  \frac{1}{n^4} \Big( \sum_{i=1}^n \sum_{j=1}^n \|X_i-X_j\| \Big)^2\\ &\leq \frac{n}{n^4} \sum_{i=1}^n \Big(\sum_{j=1}^n \|X_i-X_j\| \Big)^2 = T_{3,n}(\bX)
	\end{align*}
	and
	\begin{align*}
	T_{3,n}(\bX)&=\frac{1}{n^3} \sum_{i=1}^n \Big( \sum_{j=1}^n \|X_i-X_j\| \Big) ^2\\ &\leq \frac{n}{n^3}  \sum_{i=1}^n \sum_{j=1}^n \|X_i-X_j\|^2 = T_{1,n}(\bX).
	\end{align*}
	The second assertion in (\ref{Tineqs1and2}) follows by the triangle inequality:
	\begin{align*}
	T_{1,n}(\bX) &= \frac{1}{n^2} \sum_{i=1}^n \sum_{j=1}^n \|X_i-X_j\|^2 \\
	&= \frac{1}{n^3} \sum_{i=1}^n \sum_{j=1}^n \sum_{k=1}^n \|X_i - X_j\| \cdot \|X_i-X_k+X_k-X_j\| \\
	& \leq \frac{1}{n^3} \sum_{i=1}^n \sum_{j=1}^n \sum_{k=1}^n \|X_i - X_j\| \, \Big(\|X_i-X_k\| + \|X_k-X_j\| \Big) \\
	&= 2 \, T_{3,n}(\bX).
	\end{align*}

	The corresponding inequalities (\ref{Tineqs3and4}) for the population measures follow analogously by applying Jensen's inequality and the triangle inequality, respectively.
\end{proof}

\begin{proof}[Proof of Theorem \ref{th:dimpineq}]
	For the first assertion, we note that
	\begin{align*}
	\V_n^2(\bX) & = \big(T_{1,n}(\bX)+T_{2,n}(\bX)-2 T_{3,n}(\bX) \big) \\
	&\leq T_{2,n}(\bX) \\
	& = \frac{1}{n^4} \Big( \sum_{i=1}^n \sum_{j=1}^n \|X_i-X_j\| \Big)^2,
	\end{align*}
	where the inequality follows by (\ref{Tineqs1and2}). The second assertion follows analogously using (\ref{Tineqs3and4}).

	For establishing the last inequality, we denote the $i$-th component of $X$ and $X'$, respectively by $X^{(i)}$ and $X'^{(i)}$. Then, applying the definition of $T_1(X)$,
	\begin{align*}
	T_1(X)&= \E \|X-X'\|^2 \\
	&= \E \sum_{i=1}^p (X^{(i)}-X'^{(i)})^2 \\
	&= \sum_{i=1}^p \E \Big[(X^{(i)}-\E X^{(i)}) + (\E X^{(i)}-X'^{(i)}) \Big]^2 \\
	&= 2 \, \sum_{i=1}^p \sigma^2(X^{(i)}) = 2 \, \text{tr} (\Sigma_X).
	\end{align*}
	Applying Lemma \ref{lem:Tineq} yields
	\begin{align*}
	\V^2(X) &= T_1(X)+T_2(X)-2 T_3(X) \\
	&\leq  T_1(X)-T_3(X)  \\
	&\leq \tfrac12 \, T_1(X) \\
	&= \text{tr}(\Sigma_X).
	\end{align*}
	The proof now is complete.
\end{proof}

\begin{proof}[Proof of Proposition \ref{prop:T2T1ineq}]
	Obviously,
	\begin{equation}\label{prop:T2T1ineq1}
	1 \geq [\Cor(X,F(X))]^2 = \frac{\Cov^2(X,F(X))}{\sigma^2(X) \,\sigma^2(F(X))}.
	\end{equation}
	By \citep[equation (2.3)]{yitzhaki2003}, $\Cov(X,F(X))=\Delta(X)/4$; also, since $F(X)$ is uniformly distributed on the interval $[0,1]$ then $\Var(F(X)) = 1/12$.  By the definition of the Gini mean difference (\ref{GMD}) and by (\ref{T1T2}), $\Delta^2(X) =T_2(X)$ and $\sigma^2(X) = T_1(X)/2$.  Therefore, it follows from (\ref{prop:T2T1ineq1}) that 
	$$
	1 \ge \frac{12}{16} \, \frac{\Delta^2(X)}{{\sigma^2(X)}} = \frac{3 \,T_2(X)}{2 \, T_1(X)},
	$$
	and the proof now is complete.  
\end{proof}

\begin{proof}[Proof of Theorem \ref{th:bernoulli}]
	It is straightforward from (\ref{eq:dvartwo}) to verify that, for a Bernoulli distributed random variable $X$, 
	$\Delta(X) = 2 \, \sigma^2(X) = 2 \, T_3(X) = 2 \, p  (1-p)$.	
	%and 
	%	$T_3(X) = p (1-p)$. 	
	Hence, by (\ref{eq:dvartwo}), 
	$$
	\V^2(X) = 2 \, \sigma^2(X)+\Delta^2(X) - 2 \, T_3(X) =  4 \, p^2 (1-p)^2.
	$$	
	
	Conversely, if $X$ is a non-trivial random variable for which $\V^2(X) = \Delta^2(X)$ then the conclusion that the distribution of $X$ is concentrated on two points follows from Theorem \ref{th:repdiff}.  
\end{proof}

\begin{proof}[Proof of Proposition \ref{prop:samplelim}]
	Define
	$$
	a_{ij}=\|X_i-X_j\|, \quad b_{ij}=\|Y_i-Y_j\|, \quad \alpha_p=\E\|X - X'\|.
	$$			
	Analogous to \citep[Appendix A.1.]{szekely2013}, both $p^{-1/2} a_{12}$ and $p^{-1/2} \alpha_p$ converge, a.s., to  $\sqrt{2} \theta$ as $p \to \infty$.  Hence $a_{12}/\alpha_p \to 1$, a.s., and this yields the a.s. limits, 
	\begin{align*}
	&\lim_{p \to \infty} \frac{1}{n^2} \sum_{i,j=1}^{n} \frac{a_{ij} b_{ij}}{\alpha_p} = \frac{1}{n^2} \sum_{i,j=1}^{n} b_{ij}, \\ 
	&\lim_{p \to \infty} \frac{1}{n^4} \sum_{i,j=1}^{n} \frac{a_{ij}}{\alpha_p} \, \sum_{i,j=1}^{n} b_{ij} = \frac{n-1}{n^3} \sum_{i,j=1}^{n} b_{ij}, \\ 
	&\lim_{p \to \infty} \frac{1}{n^3} \sum_{i,j,k=1}^{n} \frac{a_{ij}  b_{ik}}{\alpha_p} = \frac{n-1}{n^3} \sum_{i,j=1}^{n} b_{ij}.	
	\end{align*}
	%	where all limits are meant a.s.
	Applying equation (\ref{eq:sampledcov}) yields assertion (\ref{eq:dcovhighdim1}). 
	
	To calculate the distance correlation, we note by \cite[Appendix A.1]{szekely2013} that 
	$$
	\lim_{p \to \infty} \frac{{\V}_n^2(\bX,\bX)}{\alpha_p^2} = \frac{n-1}{n^2},
	$$
	hence 	
	$$
	\lim_{p \to \infty} {\mathcal{R}}_n^2(\bX,\bY) = (n-1)^{-1/2} \, \frac{\frac{1}{n^2}\, \sum_{i,j=1}^{n} |Y_i-Y_j|}{{\V}_n(\bY,\bY)}.
	$$
	Applying Theorem \ref{th:dimpineq} yields
	$$
	\lim_{p \to \infty} {\mathcal{R}}_n^2(\bX,\bY) \geq  (n-1)^{-1/2},
	$$
	and this completes the proof. 
\end{proof}

\begin{proof}[Proof of Theorem \ref{th:repdiff}]
	We start by proving parts (iii) and (iv), i.e. representations (\ref{eq:dvarfour2}) and (\ref{eq:dvarfour}).
	We first derive these representations for the case in which $X$ is continuous.  In this case, we apply the Law of Total Expectation and use the mutual independence of the ranks and the order statistics \citep[Lemma 13.1]{Vandervaart2000} to obtain 
	\begin{align*}
	& W(X)\\ &=\E \Big[|X-X'| \, \big(|X-X'| - 2 \, |X-X''| \big) \Big] \\
	&= \sum_{\substack{k,k',k''=1 \\ k,k',k'' \text{are pair-} \\ \text{wise distinct}}}^3 \hspace{-12pt} \E \Big[|X-X'| \big(|X-X'| - 2 |X-X''| \big) \Big| (r_X,r_{X'},r_{X''})=(k,k',k'') \Big] \\
	& \quad \quad \quad \quad \quad \quad \quad \times \, \P \big((r_X,r_{X'},r_{X''})=(k,k',k'') \big). \\
	\end{align*}
	Using the symmetry of $X$, $X'$, and $X''$, it follows that 
	\begin{align*}
	W(X)&= \frac16 \sum_{\substack{k,k',k''=1 \\ k,k',k'' \text{are pair-} \\ \text{wise distinct}}}^3 \hspace{-12pt}  \E \Big[|X_{k:3}-X_{k':3}| \, \big(|X_{k:3}-X_{k':3}| - 2 \, |X_{k:3}-X_{k'':3}| \big) \Big] \\
	&= \frac16 \sum_{\substack{k,k',k''=1 \\ k,k',k'' \text{are pair-} \\ \text{wise distinct}}}^3  \hspace{-12pt} \E \Big[|X_{k:3}-X_{k':3}|^2 \Big] - 2 \, \E \Big[|X_{k:3}-X_{k':3}| \cdot |X_{k:3}-X_{k'':3}| \Big].
	\end{align*}
	Evaluating the first summand in the latter equation yields
	\begin{align*}
	\frac16 &\sum_{\substack{k,k',k''=1 \\ k,k',k'' \text{are pair-} \\ \text{wise distinct}}}^3  \hspace{-12pt} \E \big[|X_{k:3}-X_{k':3}|^2 \big] \\
	&= \frac13 \, \Big( \E \big[(X_{1:3}-X_{2:3})^2 \big] + \E \big[(X_{1:3}-X_{3:3})^2 \big] + \E \big[(X_{2:3}-X_{3:3})^2 \big] \Big).
	\end{align*}
	Proceeding analogously with the second summand and simplifying the outcome, we obtain	
	$$
	W(X) = - \frac43 \, \E\big[(X_{2:3}-X_{1:3})\,(X_{3:3}-X_{2:3})\big].
	$$
	This proves (\ref{eq:dvarfour}) %and (\ref{eq:dvarfour2}) 
	in the continuous case.
	
	For the case of general random variables, we now apply the method of quantile transformations. Let $U$ be uniformly distributed on the interval $[0,1]$ and let  $U$, $U'$, and $U''$ be i.i.d..  Further, let $F$ denote the cumulative distribution function of $X$.  With $F^{-1}(p) = \inf\{x: F(x) \geq p\}$ denoting the right-continuous inverse of $F$, we define $\tilde{X} = F^{-1}(\tilde{U})$, $\tilde{X'} = F^{-1}(\tilde{U'})$, and $\tilde{X''} = F^{-1} (\tilde{U''})$.  By \citep[Theorem 21.1]{Vandervaart2000}, the random variables $\tilde{X}$, $\tilde{X'}$, and $\tilde{X''}$ are i.i.d. copies of $X$ and
	{\allowdisplaybreaks
		\begin{align*}
		&W(X)\\ &=\E \Big[|\tilde{X}-\tilde{X'}| \cdot \big(|\tilde{X}-\tilde{X'}| - 2 \, |\tilde{X}-\tilde{X''}| \big) \Big] \\
		&= \sum_{\substack{k,k',k''=1 \\ k,k',k'' \text{are pair-} \\ \text{wise distinct}}}^3 \hspace{-12pt} \E \Big[|\tilde{X}-\tilde{X'}| \cdot \big(|\tilde{X}-\tilde{X'}| - 2 \, |\tilde{X}-\tilde{X''}| \big) \Big| (r_U,r_{U'},r_{U''})=(k,k',k'') \Big] \\
		& \quad \quad \quad \quad \quad \quad \quad \times \, \P \big((r_U,r_{U'},r_{U''})=(k,k',k'') \big) \\ 	
		&= \frac16 \sum_{\substack{k,k',k''=1 \\ k,k',k'' \text{are pair-} \\ \text{wise distinct}}}^3 \hspace{-12pt}  \E \Big[|X_{k:3}-X_{k':3}| \cdot \big(|X_{k:3}-X_{k':3}| - 2 \, |X_{k:3}-X_{k'':3}| \big) \Big] \\
		&= - \frac43 \, \E[(X_{2:3}-X_{1:3})\,(X_{3:3}-X_{2:3})].
		\end{align*}}
	This proves (iv), and part (iii) now follows by a combinatorial symmetry argument.

	We now prove parts (i) and (ii). For this purpose, we note that \citep{Jones2002} (see also \citep{yitzhaki2003,yitzhaki2013}) 
	\begin{equation}\label{eq:gmdrep}
	\Delta(X) = \E(D_{1:2}) = 2 \int_{-\infty}^{\infty} F(x) (1-F(x)) \dd x. 
	\end{equation}
	
	By equation (\ref{eq:gmdrep}), we obtain
	\begin{align*}
	\Delta^2(X) &= \Big[ 2 \int_{- \infty}^{\infty} F(x) \, (1-F(x)) \dd x \Big]^2  \\
	&= 4 \,\int_{- \infty}^{\infty} \int_{- \infty}^{\infty} F(x) \, [1-F(x)] \, F(y) \, [1-F(y)] \, \dd x \, \dd y \\
	&= 8 \, \operatornamewithlimits\iint\limits_{-\infty<x<y<\infty} F(x) \, [1-F(x)] \, F(y) \, [1-F(y)] \, \dd x \, \dd y. 
	\end{align*}	 	  
	Moreover, by \citep[equation (3.5)]{Jones2002}
	\begin{align*}
	& \E[(X_{2:3}-X_{1:3})\,(X_{3:3}-X_{2:3})] \\  & \quad \quad \quad \quad \quad \quad = {6} \, \operatornamewithlimits\iint\limits_{-\infty<x<y<\infty} F(x) \, [F(y)-F(x)] \, [1-F(y)] \, \dd x \, \dd y. \\
	\end{align*}
	Hence, 
	\begin{align*}
	\V^2(X) &=  \Delta^2(X)- \tfrac43 \, \E[(X_{(2)}-X_{(1)})\,(X_{(3)}-X_{(2)})] \\
	&= 8 \, \operatornamewithlimits\iint\limits_{-\infty<x<y<\infty} [F(x)]^2 \, [1-F(y)]^2 \, \dd x \, \dd y,
	\end{align*}
	which proves (ii).  
	
	Finally, (i) 
	follows from (ii) and from \citep[equation (3.4)]{Jones2002}.
\end{proof}

In proving Theorem \ref{th:disp}, we will need some preliminaries about stochastic orders of random variables 

\begin{definition}[\citep{Shaked1994}, Section 1.A.1] {\rm
		A random variable $X$ is said to be {\it stochastically smaller than} a random variable $Y$, or $X$ is {\it smaller than $Y$ in the stochastic ordering}, written $X \leq_{\text{st}} Y$, if $\P(X>u) \leq \P(Y>u)$ for all $u \in \R$.
}\end{definition}

\begin{proposition}[\citep{Shaked1994}, Section 1.A.1]
	A necessary and sufficient condition that $X \leq_{\text{st}} Y$ is that 
	\begin{equation} \label{eq:stochorder}
	\E[\phi(X)] \leq \E[\phi(Y)]
	\end{equation}	
	for all increasing functions $\phi$ for which these expectations exist.
\end{proposition}

Another important ordering of random variables is the {\it dispersive order}, $\leq_{\hskip 1pt \rm{disp}}$.

\begin{definition}[\citep{Shaked1994}, Section 2.B.1] \label{def:dispersive} {\rm
		A random variable $X$ is said to be {\it smaller than} $Y$ {\it in the dispersive ordering}, denoted by $X \leq_{\hskip 1pt \rm{disp}} Y$, if for all $0 < \alpha \leq \beta < 1$,
		\begin{equation}
		\label{eq:dispersive}
		F^{-1}(\beta) - F^{-1}(\alpha) \leq  G^{-1}(\beta) - G^{-1}(\alpha).
		\end{equation}
}\end{definition}

Bartoszewicz \citep{Bartoszewicz1986} proved the following result.

\begin{proposition}[\citep{Bartoszewicz1986}, Proposition 3] \label{th:barto}
	Let $(X_1,\ldots,X_n)$ and $(Y_1,\ldots,Y_n)$ be random samples from the random variables $X$ and $Y$, respectively, and let $D_j = X_{j+1:n}-X_{j:n}$ and $E_j = Y_{j+1:n}-Y_{j:n}$, $j=1,\ldots,n-1$ denote the corresponding sample spacings.  If $X \leq_{\hskip 1pt \rm{disp}} Y$ then $D_{j:n} \leq_{\hskip 1pt \rm{st}} E_{j:n}$ for all $j = 1,\ldots,n-1$.
\end{proposition}

\begin{proof}[Proof of Theorem \ref{th:disp}] 
	(C1) and (C2) have already been proved by Sz\'{e}kely, et al. \citep[Theorem 4]{szekely2007}. Using Definition \ref{def:dispersive}, (C3) can be reformulated as
	\begin{itemize}
		\item[(C3)] $\V(X) \leq \V(Y)$ if $X \leq_{\hskip 1pt \text{disp}} Y$.
	\end{itemize}
	Consider now the i.i.d. replicates $(X,Y)$, $(X',Y')$, $(X'',Y'')$, and $(X''',Y''')$.  Moreover, let $X_{1:4} \leq X_{2:4} \leq  X_{3:4} \leq X_{4:4}$ and $Y_{1:4}\leq Y_{2:4}\leq Y_{3:4}\leq Y_{4:4}$ denote the respective order statistics.  By Proposition \ref{th:barto}, 
	$$
	(X_{3:4}-X_{2:4}) \leq_{\hskip 1pt \rm{st}} (Y_{3:4}-Y_{2:4}). 
	$$
	Applying equation (\ref{eq:stochorder}) and the representation of the distance variance given in (\ref{eq:dvarrepcum2}) concludes the proof.
\end{proof}

\begin{proof}[Proof of Corollary \ref{cor:logconcave}]
	According to Shaked and Shantikumar \cite[Theorem 3.B.7]{Shaked1994}, a random variable $X$ satisfies the property 
	$$
	X \leq_{\hskip 1pt \rm{disp}} X+Y \ \text{for any random variable $Y$ which is independent of $X$}
	$$
	if and only if $X$ has a log-concave density.	 
	Applying this result in Theorem \ref{th:disp} concludes the proof.	
\end{proof}

\begin{proof}[Proof of Example \ref{ex:counterex1}]
	By a straightforward calculation using (\ref{eq:dvartwo}), we obtain
	\begin{align*}
	\V^2(X+Y) &= T_1(X+Y)+T_2(X+Y)-2\,T_3(X+Y) \\
	&= 
	\frac{2}{3} + \frac{4}{9} - \frac{14}{15} = \frac{8}{45}.
	\end{align*}
	However, by Theorem \ref{th:bernoulli}, $\V^2(X) = 1/4 > \V^2(X+Y)$.
\end{proof}

\begin{proof}[Proof of Example \ref{ex:nonsym}]
	By a straightforward calculation using (\ref{eq:dvartwo}), we obtain
	$$
	\V^2(X+Y) = 8 \, (p-p^2)^2 \, \big(2 \, (p-p^2)^2-6\,(p-p^2)+2 \big)
	$$
	and 
	$$
	\V^2(X-Y) = 8 \, (p-p^2)^2  \, \big(2 \, (p-p^2)^2-2\,(p-p^2)+1 \big).
	$$
	Hence, 
	$$
	\V^2(X+Y)-V^2(X-Y) = 8 \, (p-p^2)^2 \, (1-2p)^2,
	$$	
	and this difference obviously is positive for $p \neq \tfrac12$.
\end{proof}

\begin{proof}[Proof of Proposition \ref{cor:newsampledistvar}]
	Let $h: \R^4 \mapsto \R$ be the symmetric kernel defined by 
	$$
	h(X_1,\ldots,X_4) = \frac{2}{3} (X_{3:4}-X_{2:4})^2,
	$$
	where $X_{1:4} \leq X_{2:4} \leq X_{3:4} \leq X_{4:4}$ are the order statistics of $X_1,\ldots,X_4$.
	By equation (\ref{eq:dvarrepcum2}), we have $\E [h(X_1,\ldots,X_4)] < \infty$. Hence,  by Hoeffding \citep{Hoeffding1961}, 
	$$
	\widehat{\U}_n^2(\bX) = {n \choose 4}^{-1} \sum_{1\leq i_1 < i_2 < i_3 < i_4 \leq n} h(X_{i_1},\ldots,X_{i_4}) 
	$$
	is a strongly consistent estimator for $\V^2(X)$.  Using a straightforward combinatorial calculation, we obtain
	$$
	\widehat{\U}_n^2(\bX) = \frac{2}{3}{n \choose 4}^{-1} \sum_{1 \leq i < j \leq n} (i-1) \, (n-j) (X_{j:n}-X_{i:n})^2. 
	$$
	On inserting the definition of the spacings, the latter equation reduces to 
	\begin{align*}
	\widehat{\U}_n^2(\bX) &= \frac{2}{3} {n \choose 4}^{-1} \sum_{1 \leq i < j \leq n} (i-1) \, (n-j) \, (D_{i:n} + \cdots + D_{j-1:n})^2 \\
	&\equiv \frac{2}{3}{n \choose 4}^{-1} \sum_{1 \leq i < j \leq n}  (i-1) \, (n-j) \, \sum_{k,l=i}^{j-1} D_{k:n} D_{l:n}.
	\end{align*}
	Interchanging the above summations, we obtain 
	\begin{align*}
	\widehat{\U}_n^2(\bX) &= \frac{2}{3}{n \choose 4}^{-1} \sum_{k,l=1}^{n-1} D_{k:n} D_{l:n} \sum_{i=1}^{\min(k,l)} \sum_{j = \max(k,l)+1}^n (i-1) \, (n-j) \\
	%	&=\frac{2}{3}  {n \choose 4}^{-1} \sum_{k,l=1}^{n-1} D_{k:n} D_{l:n} \frac{\min(k,l) \, (\min(k,l)-1)}{2} \frac{(n-\max(k,l)) \, (n-\max(k,l)-1)}{2}\\
	&= \frac{1}{6} {n \choose 4}^{-1} \sum_{k,l=1}^{n-1} D_{k:n} D_{l:n} \min(k,l) \, \big(\min(k,l)-1\big) \\
	& \qquad\qquad\qquad\qquad\qquad \times \big(n-\max(k,l)\big) \, \big(n-\max(k,l)-1\big),
	\end{align*}
	where the latter equality follows from the fact that $\sum_{i=1}^k i = k (k+1)/2$. Since 
	$$
	\frac{1}{6} {n \choose 4}^{-1} = \frac{4}{n \, (n-1) \, (n-2) \, (n-3)},
	$$
	then we deduce that $\U_n^2(\bX) = \widehat{\U}_n^2(\bX)+o(1)$. This completes the proof. 	\end{proof}

%%%%%%%%%%%
%
%
%
%
%
%
%
%
%
\section{Remarks on the derivation of Table \ref{tab:estimators.asymptotic}}
\label{appsec:table1}

The population values and asymptotic variances of the standard deviation, the mean deviation and Gini's mean difference at the distributions considered in Table \ref{tab:estimators.asymptotic} are given by Gerstenberger and Vogel \citep[Tables 2 and 3]{gerstenberger2015}. The asymptotic variance of the maximum likelihood estimator of the scale parameter for the $t_\nu$-distribution is $(\nu + 3)/2\nu$ (see, e.g., the appendix of 
\cite{gerstenberger2018}). The population values of the distance standard deviation at the normal and the Laplace distribution are given in Appendix \ref{sec:distr}. The population value of the distance standard deviation for the $t$-distributions and its asymptotic variances at the various distributions considered can be obtained by means of numerical integration. 

For the asymptotic variance, we can employ representation (\ref{eq:asv}). However, this involves four-dimensional integrals. Generally, the numerical stability and efficiency of numerical integration quickly deteriorates as the dimension increases. We can reduce the order of integration by one by using a different representation of the asymptotic variance of the distance variance. Let
\[
\widetilde\V_n(\bX) = \widetilde{W}_n(\bX) + \widetilde{\Delta}_n^2(\bX)
\]
with
\[
\widetilde{W}_n(\bX) = \frac{1}{n(n-1)(n-2)} \sum_{\substack{1 \le i,j,k \le n \\ i\neq j, j \neq k, k \neq i}} \| X_i - X_j\|\left( \|X_i-X_j\| - 2 \|X_i- X_k\| \right)
\]
and
\[
\widetilde{\Delta}_n^2(\bX) = \frac{2}{n\,(n-1)} \sum_{1\leq i < j \leq n} \|X_i - X_j\|.
\]
The latter is the empirical Gini's mean difference, cf.~(\ref{eq:empgini}), extended to multivariate observations by replacing the absolute value with the vector norm. 
This is yet another version of the empirical distance variance. It is consistent for $\V^2(X)$ and has the same asymptotic variance as $\V_n^2(\bX)$ and $\widehat\V_n^2(\bX)$. It is not a U-statistic itself, but can be written as a function of the bivariate U-statistic
$\widetilde{B}_n(\bX) = \big(\widetilde{W}_n(\bX),\widetilde{\Delta}_n(\bX)\big)^t$ with the permutation-symmetric kernel of order three, $\tilde{h}:\R^p \times \R^p \times \R^p \to \R^2$, which maps $(x,y,z) \in \R^p \times \R^p \times \R^p$ to 
\[
%h(x,y,z) = 
\frac{1}{3}
\begin{pmatrix}
\|x\!-\!y\|(\|x\!-\!y\|\!-\!2\|x\!-\!z\|) + \|y\!-\!z\|(\|y\!-\!z\|\!-\!2\|y\!-\!x\|) \phantom{y\!-\!z \!-\!2\|y\!-\!x\|} \\
\hfill + \|z\!-\!x\|(\|z\!-\!x\|\!-\!2\|z\!-\!y\|) \\[1.0ex]
\|x-y\| + \|y-z\| + \|z-x\| \\
\end{pmatrix}.
\]
The linear part of the Hoeffding decomposition of this kernel is
\[
\widetilde{h}_1(x) = \frac{2}{3}
\begin{pmatrix}
\psi_1(x) - \psi_2(x) - 2 \, \psi_3(x) - T_1(X) + 3 \, T_3(X) \\
\psi_4(x) -T_2(X)
\end{pmatrix},
\qquad x \in \R^p,
\]
where
\[
\begin{array}{ll}
\psi_1(x) = \E \| x - Y\|^2, & \psi_2(x) = \E \|x-Y\| \|x-Z\|, \\
\psi_3(x) = \E\|x-Y\|\| Y-Z\|, & \psi_4(x) = \E\|x-Y\|,
\end{array}
\]
and $T_1(X)$, $T_2(X)$, and $T_3(X)$ are given in (\ref{T1T2}) and (\ref{T3}). Applying \cite[Theorem 7.1]{Hoeffding1948} yields that the asymptotic variance matrix of the U-statistic $\widetilde{B}_n(\bX)$ is $9 \, \E \widetilde{h}_1(\bX) \widetilde{h}_1(\bX)^t$. Denote this symmetric $2 \times 2$ matrix by $M = (m_{ij})_{i,j = 1,2}$. Its elements $m_{11}$, $m_{12}$, and $m_{22}$ are given by
\[ %begin{equation} 	\label{Mmatrix}
\begin{aligned}
m_{11} &= 4\, \E \big[ \psi_1(X) - \psi_2(X) - 2\psi_3(X) \big]^2 - 4\big(T_1(X) - 3 T_3(X)\big)^2, \\
m_{12} &= 4\, \E \big[ \psi_4(X) \big(\psi_1(X) - \psi_2(X) - 2\psi_3(X)\big) \big] - 4T_2(X)\big(T_1(X) - 3 T_3(X)\big), \\
m_{22} &= 4\, \E \psi_4^2(X) - 4 \big(T_2(X)\big)^2.
\end{aligned}
\] %\end{equation}
Applying the delta method to the function $g(x,y) = x + y^2$, we obtain %the following expression for the asymptotic variance:
\[
\gamma = m_{11}+4 \, \widetilde{\Delta}(X) \, m_{12} + 4 \,  \widetilde{\Delta}^2(X) \, m_{22},
\]
where $ \widetilde{\Delta}(X) = E\| X - X'\|$ for $X, X'$ i.i.d.\ is the population version of the multivariate Gini's mean difference. 

\section{A limit theorem for the squared distance covariance} \label{app:clt}

In this section, we state a limit theorem for the U-statistic version of the squared distance covariance, which is provided in similar form in \cite[Theorem 4.11]{Huang2017}. However,  we prove in the following that this limit theorem holds under weaker moment assumptions than known previously (cf.~\cite[Lemma 4.8]{Huang2017}). The limit theorem for independent $X$ and $Y$ (cf.~\cite[Theorem 4.12]{Huang2017}) is stated for sake of completeness. Throughout this section, we will assume that $(X,Y) \in \R^{p+q}$ are jointly distributed random vectors (with $X \in \R^p$ and $Y \in \R^q$) and for $n \in \N$, $(\bX,\bY) = ((X_1,Y_1),\ldots,(X_n,Y_n))$ are i.i.d. samples drawn from $(X,Y)$.

Moreover, for the purpose of formulating this limit theorem, we note that \cite{Huo2015} $\widehat{\Omega}_n(\bX,\bY)$ can be written as a U-statistic of order $4$ with kernel function
\begin{equation*}
\begin{aligned}
g((X_1,Y_1),\ldots,(X_4,Y_4)) &:= \frac{1}{4} \, \sum_{\substack{1 \leq i,j \leq 4 \\ i \neq j}} \|X_i - X_j \| \|Y_i - Y_j \| \\ 
& \quad - \frac{1}{4} \, \sum_{i=1}^4 \Bigg(\sum_{\substack{j=1 \\ j \neq i}}^4 \|X_i - X_j \| \sum_{\substack{j=1 \\ j \neq i}}^4 \|Y_i - Y_j \|\Bigg)  \\
& \quad + \frac{1}{24} \Bigg(\sum_{\substack{1 \leq i,j \leq 4 \\ i \neq j}} \|X_i - X_j \| \Bigg) \Bigg(\sum_{\substack{1 \leq i,j \leq 4 \\ i \neq j}} \|Y_i - Y_j \| \Bigg). \label{eq:kerneldcov} 
\end{aligned}
\end{equation*}

\begin{theorem} \label{th:asymptoticsdcov}
	Assume that $\E(\|X\|^2) < \infty$,  $\E(\|Y\|^2) < \infty$ and that $\E [g_1^2(X,Y)] > 0$, where 
	$$
	g_1(x,y) = \E [g((x,y),(X_2,Y_2),(X_3,Y_3),(X_4,Y_4))] - \V^2(X,Y).
	$$ 
	Then, for  $n \to \infty$, 
	\begin{align}
	\sqrt{n}\big(\widehat{\Omega}_n(\bX,\bY)- \V^2(X,Y)\big) \cid N(0,\,  16 \, \E [g_1^2(X,Y)]).
	\end{align}
	If $X$ and $Y$ are independent then, under existence of finite first moments (i.e., $\E(\|X\|) < \infty$ and $\E(\|Y\|) < \infty$), 
	\begin{align}
	n \, \big(\widehat{\Omega}_n^2(\bX,\bY) - \V^2(X,Y)\big) \cid \sum_{i=1}^\infty \lambda_i (Z_i^2 - 1),
	\end{align}
	where $Z_1,\ldots$ are i.i.d. standard normally distributed random variables and $\lambda_1,\ldots$ are the eigenvalues of the operator $G$ with
	$$
	G f(x_1,y_1) = \E  [6 \, g_2((x_1,y_1),(X_2,Y_2)) \, f(x_1,y_1)],
	$$
	where 
	$$
	g_2((x_1,y_1),(x_2,y_2)) = \E [g((x_1,y_1),(x_2,y_2),(X_3,Y_3),(X_4,Y_4))].
	$$ 
\end{theorem}

\begin{proof}Throughout the proof, we will denote the kernel function of the U-statistic $\widehat{\V}_n(\bX)$ by $h(\cdot)$ and the kernel function of the U-statistic $\widehat{\V}_n(\bY)$ by $\overline{h}(\cdot)$. Note that $h(\cdot)$ and $\overline{h}(\cdot)$ coincide only if $p=q$.
	
	From \cite{Huo2015}, it follows that
	$$
	g((X_1,Y_1),\ldots,(X_4,Y_4)) =	 \widehat{\Omega}_4 ((X_1,Y_1),\ldots,(X_4,Y_4)) = \frac{1}{4} \sum_{i,j=1}^4  \widetilde{A}_{ij} \widetilde{B}_{ij},
	$$
	$$ 
	h(X_1,\ldots,X_4) = \widehat{ \V}_4^2(X_1,\ldots,X_4) =	\frac{1}{4} \sum_{i,j=1}^4  \widetilde{A}_{ij}^2
	$$
	and 
	$$
	\overline{h}(Y_1,\ldots,Y_4) =	\widehat{\V}_4^2(Y_1,\ldots,Y_4) =\frac{1}{4} \sum_{i,j=1}^4  \widetilde{B}_{ij}^2,
	$$
	where $\widetilde{A}_{ij}$ is defined by
	\begin{equation}
	\widetilde{A}_{ij} = 
	\begin{cases} 
	a_{ij}-\frac{1}{2} \sum_{k=1}^n a_{kj} - \frac{1}{2} \sum_{l=1}^n a_{il} + \frac{1}{6} \sum_{k,l=1}^n a_{kl}, & i \neq j, \\ 
	0, & i=j,
	\end{cases}
	\end{equation}
	and $\widetilde{B}_{ij}$ is defined similarly.
	
	Applying twice the Cauchy-Schwarz inequality, we obtain
	\begin{align*}
	\E [g^2((X_1,Y_1),\ldots,(X_4,Y_4))] &\le \E [h(X_1,\ldots,X_4) \, \overline{h}(Y_1,\ldots,Y_4) ]  \\
	&\le \big( \E [h^2(X_1,\ldots,X_4)]  \E [\overline{h}^2(Y_1,\ldots,Y_4)] \big)^{1/2}
	< \infty,
	\end{align*}
	where the last line follows from $\E \|X\|^2<\infty$, $\E \|Y\|^2<\infty$ and  (\ref{eq:h2bounded}). Applying \citep[Theorem 7.1]{Hoeffding1948} yields 
	$$
	\sqrt{n}\big(\widehat{\Omega}_n(\bX,\bY) - \V^2(X,Y)\big) \xrightarrow{d}  \mathcal{N} \big(0, 16 \, \E [g_1^2(X)] \big).
	$$
	The limit distribution now follows  by Hoeffding \citep[Theorem 7.1]{Hoeffding1948}.
	
	For independent $X$ and $Y$,
	\begin{align*}
	&\E [g^2((X_1,Y_1),\ldots,(X_4,Y_4))] \leq \E [h(X_1,\ldots,X_4) \, \overline{h}(Y_1,\ldots,Y_4)]  \\
	&  = \E [h(X_1,\ldots,X_4)] \, \E [\overline{h}(Y_1,\ldots,Y_4)] = \V^2(X) \,  \V^2(Y)
	< \infty.
	\end{align*}
	
	Moreover, by \cite[Lemma 4.10]{Huang2017}, $g_1(x,y) = 0$ if $X$ and $Y$ are independent. The limit distribution of $\Omega_n(\bX,\bY)$ now follows by \cite[Chapter 5.5.2]{Serfling}, see also \cite[Theorem 4.12]{Huang2017}.
\end{proof}

%\newpage
\FloatBarrier

\section{Simulation results for permutation based two-sample scale tests} \label{app:permtest}

Since the small-sample distribution of $|\widehat{T}_\V|$ is complicated, permutation tests provide a practical approach to obtaining its small-sample critical values.  To derive critical values of the test, we implement in the simulations below the following permutation sampling scheme: 
%The maximal number of permutations is $1000$, i.e., 
For i.i.d., mutually independent samples $\bX_n=(X_1,\ldots,X_n)$ and $\bY\!_m= (Y_1,\ldots,Y_m)$, we draw without replacement a random sample of size $n$ from the $n+m$ data points, label this sample $\tilde\bX_n$, and label the remaining $m$ data points $\tilde\bY\!_m$, and compute the test statistic from the two new samples. This is repeated $1,000$ times, and the proportion of times that the permutation test statistic exceeds the original test statistic yields the $p$-value of the test. %For $\binom{n+m}{n} \le 1000$, instead of random sampling, all possible permutations are considered.

We remark that, different from the asymptotic test described in Section \ref{subsec:hyptesting} in the main paper, the permutation based test can only be applied under the assumption that the distributions corresponding to the samples $\bX_n$ and $\bY\!_m$ belong to the same location-scale family and share the same location. 

We consider the same simulation scenarios as for Tables \ref{tab:tests.size} and \ref{tab:tests.power} in the main paper. Specifically, we consider the Laplace distribution, normal distribution, normal scale mixture distribution $NM(3,0.1)$, and the $t_\nu$-distributions with $\nu = 3$ and $\nu = 5$.
The sample sizes $n, m$ range from $n + m = 30$ to $n + m = 2,000$. Table \ref{tab:tests.sizeperm} (test size) contains results for the null hypothesis $\lambda = 1$ and Table \ref{tab:tests.powerperm} (test power) gives results for the sample-size-dependent alternative with
\begin{equation} \nonumber	%\label{lambdanm}
\lambda_{n, m} = 1 + 3\sqrt{\frac{n+m}{n \, m}}. 
\end{equation}

In Table \ref{tab:tests.sizeperm}, we observe that the actual rejection frequencies of all tests are close to $5.0$, as one can expect from permutation-based methods.

Table \ref{tab:tests.powerperm} illustrates that the permutation-based distance standard deviation test shows a considerably better performance than the asymptotic test for the small-sample setting $(n,m) = (15,15)$.  With increasing sample sizes, the advantages of the permutation-based approach get smaller; for sample sizes greater than $n+m = 500$, the two approaches perform almost equally.

Comparing the performance of the tests based on $|\widehat{T}_\V|$, $|\widehat{T}_\Delta|$ and $|\widehat{T}_\sigma|$ yields very similar results as obtained from comparing the corresponding asymptotic tests in the main paper. 

\begin{table}[t!]
	\captionsetup{width=\textwidth}
	\caption{\small \emph{Test size.}  Empirical rejection frequencies (\%) under the null hypothesis $\lambda=1$ of permutation-based two-sample scale tests (based on the distance standard deviation $\widehat\V_n$, the standard deviation $\widehat\sigma_n$, Gini's mean difference $\widehat\Delta_n$, and the $F$-test) at the 5\% significance level. Results are based on $10,000$ replications.
		\label{tab:tests.sizeperm}}
	\renewcommand{\arraystretch}{1.05}
	%\small
	\begin{center}
		\begin{tabular}{cc|cccccc} 
			& $n$ 				& 15 & 50 & 120 	& 250 & 600 & 1,000 \\
			& $m$ & 15  & 50 	& 40 		& 250 		& 200 		& 1,000 			 \\[1.1ex]
			Distribution & Test & \multicolumn{5}{c}{Rejection frequencies (\%)} \\[1.2ex]
			\hline {\phantom{\Big[}}
			$L(0,1)$ 	& $\widehat\V_n$  		& 5.5 &  4.6 &   5.0 &    5.0  &   4.9  &     4.7 \\
			& $\widehat\sigma_n$ 	& 5.4 &  4.8  &  5.2  &   4.5  &   5.0  &     5.1\\
			& $\widehat\Delta_n$ 	& 5.3 &   4.8  &  5.2  &   4.8   &  4.9   &   4.8\\
			$N(0,1)$ 	& $\widehat\V_n$ 			& 5.1 &  5.0 &   5.1 &    4.9  &   4.7  &     4.9\\
			& $\widehat\sigma_n$ 	& 4.9 &  4.9 &   5.2 &    5.2 &    5.2  &     4.8\\
			& $\widehat\Delta_n$	& 5.0 &   4.8 &   5.2  &    4.9 &  5.1 &      4.8\\
			$NM(3,0.1)$ & $\widehat\V_n$ 		& 4.7 &  5.1 &   5.1  &   5.2  &   5.3  &     5.2\\
			& $\widehat\sigma_n$ 	& 4.7 &  5.2  &  5.0 &    5.0 &    5.1  &     5.1\\
			& $\widehat\Delta_n$ 	& 4.5  &   5.3  &  5.0 &    4.8  &   4.9  &    5.1\\
			$t_3$ 		& $\widehat\V_n$ 			& 5.5 &  4.9  &  4.8  &   4.6  &   4.7  &     4.8\\
			& $\widehat\sigma_n$ 	& 5.2  & 5.0 &   5.2  &   5.1  &   5.0   &   4.9 \\
			& $\widehat\Delta_n$ 	& 5.3 &   5.1 &   5.1 &    4.7  &   5.1  &     4.7\\
			$t_5$ 		& $\widehat\V_n$ 	&		5.1  & 5.2   & 5.0  &   5.1 &    5.3   &    4.9\\
			& $\widehat\sigma_n$ 	& 5.2 &  5.2   & 5.0 &    4.8 &    5.1    &   4.9\\
			& $\widehat\Delta_n$ 	& 5.2  &  5.3  &  4.9  &   5.1  &   5.3  &     5.1\\
			\hline
		\end{tabular}
	\end{center}
\end{table}

\begin{table}[ht!]
	\captionsetup{width=\textwidth}
	\caption{\small \emph{Test power.} Empirical rejection frequencies (\%) under the alternative $\lambda_{n,m} = 1 + 3 \sqrt{(n+m)/n/m}$ of permutation-based two-sample scale tests (tests based on the distance standard deviation $\widehat\V_n$, the standard deviation $\widehat\sigma_n$, Gini's mean difference $\widehat\Delta_n$, and the $F$-test) at the 5\% significance level. Results are based on $10,000$ replications.		
		\label{tab:tests.powerperm}}
	\renewcommand{\arraystretch}{1.05}
	%\small
	\begin{center}
		\begin{tabular}{cc|cccccc} 
			& $n$ 				& 15 & 50 & 120 	& 250 & 600 & 1,000 \\
			& $m$ & 15  & 50 	& 40 		& 250 		& 200 		& 1,000 			 \\[1.1ex]
			Distribution & Test & \multicolumn{5}{c}{Rejection frequencies (\%)} \\[1.2ex]
			\hline  {\phantom{\Big[}}
			$L(0,1)$ 	& $\widehat\V_n$  		& 41.0 & 59.6 &  65.1  &  73.4  &  75.4  &   78.7 \\
			& $\widehat\sigma_n$ 	& 41.9 & 57.0 &  61.9 &   66.6  &  69.6  &    71.5\\
			& $\widehat\Delta_n$ 	& 43.6 & 60.7 &  68.0 &   73.6  &  77.0  &   79.1\\
			$N(0,1)$ 	& $\widehat\V_n$ 			& 58.2 & 80.5 &  84.5  &  90.8 &   92.3  &    94.3\\
			& $\widehat\sigma_n$ 	& 67.9 &  88.3 &  91.5  &  96.0 &   96.5  &   97.8\\
			& $\widehat\Delta_n$	& 67.4 & 87.8 &  91.4  &  95.5 &   96.4   &  97.6\\
			$NM(3,0.1)$ & $\widehat\V_n$ 		& 50.3 &  70.6 &  75.1 &   82.5  &  84.5 &     87.5\\
			& $\widehat\sigma_n$ 	& 50.0 & 51.9 &  53.3 &   52.2   & 54.4  &    55.1\\
			& $\widehat\Delta_n$ 	& 51.2 & 61.1  & 67.1  &  69.6   & 73.3 &     75.8\\
			$t_3$ 		& $\widehat\V_n$ 			&44.3&  61.7 &  65.3 &   72.8  &  76.1  &    78.9\\
			& $\widehat\sigma_n$ 	& 42.5 & 46.5 &  44.0 &   38.4 &   34.5  &    29.8 \\
			& $\widehat\Delta_n$ 	& 44.3 & 53.5  & 56.2  &  57.2 &   58.1   &   58.1\\
			$t_5$ 		& $\widehat\V_n$ 	&	50.5 & 70.4  & 76.1  &  83.8 &   84.5   &   87.9\\
			& $\widehat\sigma_n$ 	& 52.1 &  63.9 &  66.4 &   67.7   & 66.1   &   66.4\\
			& $\widehat\Delta_n$ 	& 53.2 &  68.9 &  74.9  &  79.2 &   80.4   &   83.4\\
			\hline
		\end{tabular}
	\end{center}
\end{table}	

\FloatBarrier

\section{Further properties of the distance standard deviation in one dimension} \label{app:furtherprop}

\paragraph{Consequences of representation (\ref{eq:dvarfour2})}

%\medskip

In the continuous case with finite second moment, equation (\ref{eq:dvarfour2}) is equivalent to the identity, 
\begin{equation}{\label{eq:gmdvar}}
\E ( |X-X'| \cdot |X''-X'| ) = \sigma^2(X) + 4 \, J(X),
\end{equation}
where 
$$
J(X) = \int_{x= -\infty}^\infty \int_{y =- \infty}^{x} \int_{z=x}^\infty (x-y) \, (z-x) f(z) \, f(y) f(x) \dd z \dd y \dd x	.
$$
Formula (\ref{eq:gmdvar}) is essentially the key result in the classical paper by Lomnicki \citep{lomnicki1952}, who also obtained a simple expression for the variance of the empirical Gini's mean difference $\widehat{\Delta}_n$.

Indeed, it was shown in \citep{lomnicki1952} that, for continuous random variables $X$ and an i.i.d. sample $\bX = (X_1,\ldots,X_n)$ drawn from $X$,
%the variance of the estimator is given by
\begin{multline} \label{eq:gmdvartwo}
\Var\big(\widehat{\Delta}_n(\bX) \big) \\
= \frac{1}{n \, (n-1)} \big[4 \,(n-1) \,  \sigma^2(X)+ 16 \, (n-2) J(X)- 2 \, (2 n-3) \Delta^2(X) \big].
\end{multline}

We note several consequences of Theorem \ref{th:repdiff} and equation (\ref{eq:gmdvartwo}).  First, Theorem \ref{th:repdiff} implies that the decomposition (\ref{eq:gmdvartwo}) holds in an analogous way for the non-continuous case. Second, Theorem \ref{th:repdiff} implies
$$
J(X) = \tfrac18 (\Delta^2(X) - \V^2 (X) ).
$$
Inserting this expression into (\ref{eq:gmdvartwo}) yields that the variance of $\widehat{\Delta}_n$ can be established by calculating the corresponding distance variance and {\it vice versa}. In particular, solving the integral $(\ref{eq:dvar})$ for some random variable $X$ gives us the corresponding variance of the empirical Gini's mean difference. These considerations imply that the asymptotic variance $$\mbox{ASV}(\widehat{\Delta}_n(\bX)) = \lim_{n \to \infty} n \, \Var\big(\widehat{\Delta}_n(\bX) \big) $$ can be expressed alternatively as
\begin{equation} {\label{eq:ASV}}
\mbox{ASV}(\widehat{\Delta}_n(\bX)) = 4 \, \sigma^2(X) - 2 \, \V^2 (X) - 2 \, \Delta^2(X).
\end{equation}
Using equation (\ref{eq:ASV}), novel expressions for the asymptotic variance of $\widehat{\Delta}_n$ for the gamma, the negative binomial and the Poisson distribution can be established from the results for the distance variance in \cite{Dueck2016}; similarly novel expressions for the distance variance of the uniform, the Laplace, the Pareto and the exponential distribution can be obtained from the results for the asymptotic variance of $\widehat{\Delta}_n$ in \cite{Zenga2004,gerstenberger2015}.

%By equations (\ref{eq:gmdvartwo}) and (\ref{eq:ASV}), we can also derive expressions for the variance and asymptotic variance of the sample Gini mean difference for a wide class of distributions, as stated in Theorem \ref{th:distri}.  To the best of our knowledge, those expressions are novel for the Gamma, Poisson, and negative binomial distributions.

\paragraph{Consequences of representation (\ref{eq:dvarfour})}

Representation (\ref{eq:dvarfour}) enables us to establish inequalities for the distance standard deviation of random variables with log-concave and log-convex densities. Clearly 
$$
\E[(X_{2:3}-X_{1:3})\,(X_{3:3}-X_{2:3})] \geq \E[X_{2:3}-X_{1:3}] \, \E[X_{3:3}-X_{2:3}],
$$
whenever $\mbox{Cov} (X_{2:3}-X_{1:3},\,X_{3:3}-X_{2:3}) \geq 0$.  Moreover, the reverse inequality holds if the respective covariance is smaller than or equal to $0$.  The dependence structure of spacings has been studied by Yao, et al. \cite{yao2008}; they showed that if a random variable $X$ is supported on $(a,\infty)$, where $a \in \R \cup \{-\infty\}$, and has a log-convex density then the sequence of spacings $(D_{1},\ldots,D_{n-1})$ of $X$ is {\it multivariate totally positive of order 2} ($MTP_2$). They also proved that if $X$ has a log-concave density then the sequence of spacings of $X$ is {\it multivariate reverse rule of order 2} ($MRR_2$). These relationships lead to the following result.   

\begin{proposition}{\label{prop:vg12delta}}
	Let $a \in \R \cup \{-\infty\}$ and $b \in \R \cup \{\infty\}$. If $X$ is a random variable with log-convex density and support $(a,\infty)$, then 
	$$
	\V^2(X) \leq \Delta^2(X)- \tfrac43 \, \E[X_{2:3}-X_{1:3}] \, \E[X_{3:3}-X_{2:3}].
	$$	
	Moreover if $X$ is a random variable with log-concave density and support $(a,b)$, then 
	$$
	\V^2(X) \geq \Delta^2(X)- \tfrac43 \, \E[X_{2:3}-X_{1:3}] \, \E[X_{3:3}-X_{2:3}]
	$$
	and, consequently, 
	$$
	\V^2(X) \geq \tfrac{1}{4} \Delta^2(X).
	$$		
\end{proposition}
\begin{proof}
	For the first inequality, we apply \cite[Theorem 3.1, Remark 3.3]{yao2008} to show that $(X_{2:3}-X_{1:3}, \, X_{3:3}-X_{2:3})$ is $MTP_2$ if $X$ has a log-concave density with support $(a,\infty)$. By \cite[Equation (1.7)]{karlin1980a}, this induces that
	$$
	\mbox{Cov} (X_{2:3}-X_{1:3},\,X_{3:3}-X_{2:3}) \geq 0,
	$$
	which completes the proof.
	The second inequality follows analogously by applying \cite[Theorem 3.1, Remark 3.3]{yao2008} and \cite[Lemma 2.1]{karlin1980b}. The third inequality follows by the second inequality using 
	$$
	\E[X_{2:3}-X_{1:3}] \, \E[X_{3:3}-X_{2:3}] \leq \tfrac14 \E^2[X_{3:3}-X_{1:3}] = \tfrac{9}{16} \Delta^2(X).
	$$		
\end{proof}

%	\medskip

Combining Proposition \ref{prop:vg12delta} and Corollary \ref{cor:vargmdineq}, we obtain for random variables $X$ with log-concave densities, that
$$
\frac12  \leq \frac{\V(X)}{\Delta(X)} \leq 1.
$$
Moreover, combining Proposition \ref{prop:vg12delta} and equation (\ref{eq:ASV}) yields an upper bound for the asymptotic variance of Gini's mean difference for random variables $X$ with log-concave densities:
$$
ASV(\widehat{\Delta}_n(\bX)) \leq 4 \sigma^2(X) - \tfrac52 \Delta^2(X).
$$

\section{Expressions for the distance variance of some well-known distributions}
\label{sec:distr}

Exploiting the different representations of the distance variance derived in Theorem \ref{th:repdiff}, we can state the distance variance of many well-known distributions. In the following, we tabulate these distance variances for future reference. We use the standard notation ${}_1F_1$ and ${}_2F_1$ for the classical confluent and Gaussian hypergeometric functions.  

\begin{theorem} 
	\label{th:distri}
	\begin{enumerate}
		\item Let $X$ be Bernoulli distributed with parameter $p$. Then 
		$
		\V^2(X) = 4\,p^2\,(1-p)^2.
		$
		\item Let $X$ be normally distributed with mean $\mu$ and variance $\sigma^2$. Then 
		$$
		\V^2(X) = 	4\Big( \frac{1-\sqrt{3}}{\pi} + \frac{1}{3} \Big) \sigma^2.
		$$	
		\item Let $X$ be uniformly distributed on the interval $[a,b]$. Then 
		$
		\V^2(X) = 	2(b-a)^2/45.
		$
		\item Let $X$ be Laplace-distributed with density function, $f_X(x)=(2 \alpha)^{-1}$ \\ $\exp(-|x-\mu|/\alpha)$, $x \in \R$, $\alpha > 0$, $\mu \in \R$. Then
		$
		\V^2(X) = 7\alpha^2/12.
		$
		\item Let $X$ be Pareto-distributed with parameters $\alpha > 1$ and $\lambda > 0$, and density function $f_X(x) = \alpha \lambda^\alpha x^{-(\alpha+1)}$, $x \geq \lambda$.  Then, 
		$$
		\V^2(X) = \frac{4 \alpha^2 \lambda^2}{(\alpha-1) \, (2 \alpha-1)^2 \, (3 \alpha -2)}.
		$$						
		\item Let $X$ be exponentially distributed with parameter $\lambda > 0$ and density function $f_X(x)= \lambda \exp(-\lambda \, x)$, $x \geq 0$.  Then, 
		$
		\V^2(X) = (3 \lambda^2)^{-1}.
		$		
		\item Let $X$ be Gamma-distributed with shape parameter $\alpha > 0$ and scale parameter $1$. Then
		$$
		\V^2(X) = 2^{2(2-2\, \alpha)}  \sum_{j, k=1}^\infty A_{j,k}^2(\alpha),
		$$
		where 
		\begin{align*}
		A_{j,k}&(\alpha) = 2^{-j-k} \, \left(\frac{(\alpha)_j \, (\alpha)_k}{j! \, k!} \right)^{1/2} \\
		& \times \frac{\Gamma(2 \alpha+j+k-1)}{\Gamma(\alpha+j) \, \Gamma(\alpha+k)}\ {}_2F_1\left(-j-k+2,1-\alpha-j;2-2 \alpha-j-k;2 \right). 
		\end{align*}
		\item Let $X$ be Poisson-distributed with parameter $\lambda > 0$. Then
		$$
		\V^2(X) = \sum_{j, k=1}^\infty \frac{4^{j+k-1} }{j! \, k!} \, \lambda^{j+k} \, A_{jk}^2, 
		$$
		where 
		\begin{align*}
		A_{jk} = \frac{1}{(j-1)!} \sum_{l=0}^{\lfloor (j-k)/2\rfloor} \binom{j-k}{2l} (-1)^l (\tfrac12)_l \, (\tfrac12)_{j-l-1} \ {}_1F_1(j-l-\tfrac12;j;-4\lambda).
		\end{align*}
		\item Let $X$ have a negative binomial distribution with parameters $c$ and $\beta$. Then
		$$
		\V^2(X) = (1-c)^{4\beta} \sum_{j,k=1}^\infty \frac{(\beta)_j \, (\beta)_k}{j! \, k!} (1+c^2)^{-2\beta-2j} 2^{2k}  c^{j+k} A_{jk}^2,
		$$
		where 
		\begin{align*}
		A_{jk} &= \sum_{l_1,l_2=0}^{j-k} \binom{j-k}{l_1} \binom{j-k}{l_2} (-c)^{l_1} (-1)^{l_2} (|l_1-l_2|)! \sum_{l=0}^\infty \frac{(\beta+j)-l}{l!} \bigg(\frac{2c}{1+c^2}\bigg)^l \\
		& \qquad \times \sum_{m=0}^{|l_1-l_2|} (-2)^m \frac{(m)_{|l_1-l_2|}}{(|l_1-l_2|-m)! \, (2m)!} \frac{2^{k+m-1} \, (\tfrac12)_{k+m-1}}{(k+m-1)!} \\ & \qquad \qquad \times {}_2F_1(-l,k+m-\tfrac12;k+m;2).
		\end{align*}
		\item Let $X=(X_1,\ldots,X_p)$ be a multivariate normally distributed random vector with mean $\mu=(\mu_1,\ldots,\mu_p)$ and identity covariance matrix $I_p$. Then
		$$
		\V^2(X) =    4 \pi \,  \frac{c_{p-1}^2}{c_{p}^2} \, \left[
		\frac{\Gamma(\tfrac12 p) \, \Gamma(\tfrac12 p + 1)}
		{\big[\Gamma\big(\tfrac12(p+1)\big)\big]^2} 
		- 2 \ {}_2F_1 \! \left(-\tfrac12,-\tfrac12;\tfrac12 p;\tfrac14\right) + 1\right].
		$$
	\end{enumerate}
\end{theorem}

\begin{proof}[Proof of Theorem \ref{th:distri}]
	1.  See Theorem \ref{th:bernoulli}.
	
	2. See the proof of Theorem 7 in \citep{szekely2007} or \citep[p. 14]{Dueck2014}.
	
	3. and 4.  These follow directly from Theorem \ref{th:repdiff} and the results in Table 3 in \citep{gerstenberger2015}.
	
	5. and 6.  These results follow directly from the representation (\ref{eq:dvartwo}) and \citep[equations (4.2) and (4.4)]{Zenga2004}.
	
	7., 8., and 9.  See \citep[Propositions 5.6, 5.7, and 5.8]{Dueck2016}.
	
	10. See \citep[Corollary 3.3]{Dueck2014}.
\end{proof}

\end{document}